\documentclass [12pt]{article}
\usepackage{amsmath,amssymb,amscd}
\usepackage{hyperref}
\usepackage{graphicx}

\begin{document}
\def\l{\lambda}
\def\tilde{\widetilde}
\def\m{\mu}
\def\a{\alpha}
\def\b{\beta}
\def\GCD{{\rm GCD }}
\def\gcd{{\rm GCD }}
\def\g{\gamma}
\def\G{\Gamma}
\def\d{\delta}
\def\e{\epsilon}
\def\o{\omega}
\def\O{\Omega}
\def\v{\varphi}
\def\t{\theta}
\def\r{\rho}
\def\bs{$\blacksquare$}
\def\bp{\begin{proposition}}
\def\ep{\end{proposition}}
\def\bt{\begin{theo}}
\def\et{\end{theo}}
\def\be{\begin{equation}}
\def\ee{\end{equation}}
\def\bl{\begin{lemma}}
\def\el{\end{lemma}}
\def\bc{\begin{corollary}}
\def\ec{\end{corollary}}
\def\pr{\noindent{\bf Proof: }}
\def\note{\noindent{\bf Note. }}
\def\bd{\begin{definition}}
\def\ed{\end{definition}}
\def\C{{\mathbb C}}
\def\P{{\mathbb P}}
\def\Z{{\mathbb Z}}
\def\d{{\rm d}}
\def\deg{{\rm deg\,}}
\def\deg{{\rm deg\,}}
\def\arg{{\rm arg\,}}
\def\min{{\rm min\,}}
\def\max{{\rm max\,}}
 \def\d{{\rm d}}
\def\deg{{\rm deg\,}}
\def\deg{{\rm deg\,}}
\def\arg{{\rm arg\,}}
\def\min{{\rm min\,}}
\def\max{{\rm max\,}}

\newtheorem{theo}{Theorem}[section]
\newtheorem{lemma}{Lemma}[section]
\newtheorem{definition}{Definition}[section]
\newtheorem{corollary}{Corollary}[section]
\newtheorem{proposition}{Proposition}[section]

\begin{titlepage}
\begin{center}

\topskip 2mm

{\LARGE{\bf {}Algebraic Geometry of the Center-Focus problem

\medskip

for Abel Differential Equation}}

\vskip 4mm

{\large {\bf M. Briskin $^{*}$,}}
\hspace {1 mm}
{\large {\bf F. Pakovich $^{**}$,}}
\hspace {1 mm}
{\large {\bf Y. Yomdin $^{***}$}}

\vspace{2 mm}
\end{center}

{$^{*}$ Jerusalem College of Engineering, Ramat Bet Hakerem, P.O.B. 3566,
Jerusalem 91035, Israel},
{$^{**}$ Department of Mathematics, Ben-Gurion University of the Negev,
Beer-Sheva 84105, Israel},
{$^{***}$ Department of Mathematics, The Weizmann Institute of
Science, Rehovot 76100, Israel}

\vspace{1 mm}

\noindent e-mail: briskin@jce.ac.il, \ pakovich@math.bgu.ac.il,
 \ yosef.yomdin@weizmann.ac.il

\vspace{1 mm}
\begin{center}

{ \bf Abstract}
\end{center}
{\small{

The Abel differential equation $y'=p(x)y^3 + q(x) y^2$ with polynomial coefficients
$p,q$  is said to have a center on $[a,b]$ if all its solutions, with the
initial value $y(a)$ small enough, satisfy the condition $y(a)=y(b)$.
The problem of giving conditions on $(p,q,a,b)$ implying a center for the Abel equation
is analogous to the classical Poincar\'e Center-Focus problem for plane vector fields.

Center conditions are provided by an infinite system of ``Center Equations''. During the last two decades, an important
new information on these equations has been obtained via a detailed analysis of two related
structures: Composition Algebra and Moment Equations (first order approximation of the
Center ones). Recently, one of the basic open questions in this direction - the ``Polynomial
moments problem'' - has been completely settled in \cite{mp1,pak}.

In this paper we present a progress in the following two main directions: First, we translate
the results of \cite{mp1,pak} into the language of Algebraic Geometry of the Center
Equations. Applying these new tools, we show that the center conditions can be described in terms
of Composition Algebra, up to a ``small'' correction. In particular, we significantly
extend the results of \cite{broy}. Second, applying these tools in combination with explicit computations,
we start in this paper the study of the ``second Melnikov coefficients'' (second order approximation of
the Center Equations) showing that in many cases vanishing of the moments and of these coefficients
is sufficient in order to completely characterize centers.

\vspace{1 mm}
\begin{center}
------------------------------------------------
\vspace{1 mm}
\end{center}
This research was supported by the ISF, Grants No. 639/09 and No. 779/13, and by the
Minerva Foundation.

}}

\end{titlepage}
\newpage


\section{Introduction}
\setcounter{equation}{0}

In this paper we consider the Abel differential equation
\be \label{Abel}
y'=p(x)y^3 + q(x) y^2
\ee
with polynomial coefficients $p,q$, on a complex segment $[a,b]$. A solution $y(x)$ of (\ref{Abel}) is called
``closed'' on $[a,b]$ if $y(a)=y(b)$ for the initial element of $y(x)$ around $a$ analytically continued to $b$ along $[a,b]$.
Equation (\ref{Abel}) is said to have a center on $[a,b]$ if any its solution $y(x)$, with the initial value $y(a)$ small
enough, is closed on $[a,b]$. For $p,q$ polynomials this property depends only on the endpoints $a,b\in {\mathbb C}$, but not
on the continuation path.

Below we shall denote by $P,Q$ the primitives $P(x)=\int_a^x p(\tau)d\tau$ and $Q(x)=\int_a^x q(\tau)d\tau$.

\smallskip

The Center-Focus problem for the polynomial Abel equation is to give an explicit, in terms of the coefficients of $p$ and $q$,
necessary and sufficient condition on $p,q,a,b$ for (\ref{Abel}) to have a center on $[a,b]$. The Smale-Pugh problem is to bound
the number of isolated closed solutions of (\ref{Abel}). While we restrict ourselves to the polynomial case only, there are
other important settings of these problems, in particular, with $p,q$ trigonometric polynomials, piecewise-linear, or even
discontinuous piecewise constant functions (compare \cite{Alv.Bra.Chr,al,al3,bru,Bru1,cgm1,cgm2,cgm3,ggl,ggs1,ggs2}). The relation 
of the above problems to the classical Hilbert 16-th and Poincar\'e Center-Focus problems for plane vector fields is well known 
(see, eg. \cite{broy,chr,Ily,lpr}).

Algebraic Geometry enters the above problems from the very beginning:
it is well known that center conditions are given by an infinite system of
polynomial equations in the coefficients of $p,q$, expressed as certain
iterated integrals of $p,q$ (``Center Equations''; see Section 3 below).
The structure of the ideal generated by these equations in an appropriate
ring (called the Bautin ideal), specifically, the number of its generators, determines local bifurcations of the closed solutions
as $p,q$ vary.

One of the main difficulties in the Center-Focus and the Smale-Pugh  problems
is that a general algebraic-geometric analysis of the system of Center
Equations is very difficult because of their complexity and absence of
apparent general patterns.

In recent years the following two important
algebraic-analytic structures, deeply related to the Center Equations
for (\ref{Abel}), have been discovered: Composition Algebra of polynomials
and generalized polynomial moments of the form
$m_k=\int_a^b P^k(x)q(x)dx$ (the last one is a special
case of iterated integrals). The use of these structures provides important tools
for investigation of the Center-Focus problem for the Abel equation
(see \cite{Alv.Bra.Chr}-\cite{cgm3}, \cite{gssy,ggl,ggs1,ggs2,ny} and references therein).
In particular, it was shown in \cite{broy} that Center Equations are well
approximated by the Moment Equations $m_k=\int_a^b P^k(x)q(x)dx=0$, and
in fact coincide with them ``at infinity''.  Moment Equations, in turn,
impose (in many cases) strong restrictions on $P$ and $Q$, considered as elements
of the Composition Algebra of polynomials (see Section \ref{Mom.Comp} below). Notice 
that usually {\it linear} Moment Equations
$m_k=0$ are considered, where $P$ is fixed while $Q$ is the unknown.
However, consideration of Center Equations at infinity in \cite{broy,gssy} and in
the present paper leads to a {\it non-linear} setting where $Q$ is fixed, while
the equations have to be solved with respect to the unknown $P$.

The following Composition Condition  imposed on $P$ and $Q$ plays a
central role in the study of the Moment and Center Equations (see the references above):
there exist polynomials $\tilde P, \ \tilde Q$ and $W$ with
$W(a)=W(b)$ such that $$P(x)=\tilde P(W(x)), \ Q(x)=\tilde Q(W(x)).$$
Being a kind of ``integrability condition'', Composition Condition implies vanishing
of Center and Moments Equations as well as of all the iterated integrals entering
the Center Equations. It is the only known to us {\it sufficient} center condition
for the polynomial Abel equation. Using the interrelation between Center and Moment
Equations at infinity, and the Composition Condition, a rather accurate description
of the affine Center Set for the polynomial Abel equation has been given in \cite{broy}. Very
recently important results relating Center and Composition Conditions for trigonometric
and polynomial Abel equations have been obtained in \cite{bpy,cgm1,cgm2,cgm3,ggl,ggs1,ggs2}.

\smallskip

These results, as well as some further examples and partial results (see \cite{al2,bby,broy,bru,Bru1,Bru.Yom,cgm2,gssy}
for the most recent contributions) seem to support the following ``Composition conjecture'':

\smallskip

\noindent{\bf Conjecture 1.} {\it The Center and Composition Sets for any
polynomial Abel equation coincide}.

\smallskip

This conjecture was originally suggested in \cite{bfy} (Conjecture 1.6), together with its extended versions (\cite{bfy}, 
Conjectures 1.7 and 1.8) which all remain open. Similar conjecture is known to be false for $p,q$ trigonometric polynomials 
and $a,b\in {\mathbb R}$ (see \cite{al}). However, besides various special cases of polynomial Abel equations, described in 
papers mentioned above, as well as in \cite{ggl,ly} and in other publications, an equivalence of the Center and (the appropriate)
Composition Conditions holds, for example, for piecewise-constant $p$ and $q$ of a certain special
form (``rectangular paths'', \cite{Bru1}, see also \cite{al3}). As it was shown in \cite{Bru1}, for ``rectangular paths'' the
equivalence of the Center and Composition conditions follows from a highly non-trivial result of \cite{Coh}, stating (roughly) 
that the group of transformations of ${\mathbb R}$ generated by translations and positive rational powers is free.

\smallskip

Part of the methods developed in the present paper can be applied to arbitrary coefficients $p,q$ of Abel equation
(\ref{Abel}). This certainly concerns all the constructions in Section \ref{Cen.Eq} below. In particular, we can apply our
methods to $p,q$ - trigonometric polynomials, Laurent polynomials, or rational functions. The problem is that in the
case of rational $p,q$ the consequences of the moments vanishing are much weaker than in the polynomial case, while
the presentation is technically much more involved (see \cite{Alv.Bra.Chr,pak1} and references therein). The same is true for 
the description of the Composition algebra of rational functions, which turns out to be significantly more complicated than
for polynomials (compare \cite{Alv.Bra.Chr,cgm1,cgm2,pak2,pak3}). So in the present paper we restrict ourselves to the 
polynomial case only. We plan to present our results for rational and trigonometric cases separately.

\smallskip

Now, in \cite{mp1,pak} essentially a complete description of the {\it polynomial} moments vanishing has been achieved,
as well as of the relevant {\it polynomial} Composition Algebra. In particular, explicit necessary and sufficient
conditions for vanishing of all the moments $m_k$ have been given there, in terms of certain relations between $P,Q,a,b$
in the Composition Algebra of polynomials (see Section \ref{Mom.Comp} below).

\medskip

Accordingly, one of the main goals of the present paper is to give an algebraic-geometric interpretation of the results
of \cite{pp,mp1,pak} in the context of the Center-Focus problem for the polynomial Abel equation, and to apply these
results to the study of Center Conditions. Here we heavily use the fact, found in \cite{broy}, that the Moment Equations
are the restrictions, in a proper ``projective setting'' of the Center Equations to the infinite hyperplane. On this base
we obtain new information on the affine Center Conditions, significantly extending the results of \cite{broy}.

Another main goal of this paper is to start the investigation of the ``second Melnikov
coefficients'', which form the second set of the Center Equations ``at infinity''. We show that
in many important cases {\it vanishing of the moments and of the second Melnikov coefficients implies
composition}, and so it is sufficient in order to completely characterize centers.

\smallskip

The authors would like to thank the referee for detailed suggestions, significantly improving the presentation.

\subsection{Statement of the main results} \label{st.main}

A general form of the results in this paper is the following: as it was explained above, Composition Set is always a subset
of the Center Set. We show that the Composition Condition is indeed a good approximation to the Center Condition, showing 
that the dimension of the (possibly existing) non-composition components in the Center Set is small. In various circumstances 
we provide an upper bound for the dimension of these possible non-composition components, which is significantly smaller
than the dimension of the Composition Center strata. In many cases this bound is zero, so the Center Set coincides with the 
Composition Set up to a finite number of points. The following theorems summarize our main new results on the center configurations 
for polynomial Abel equation (\ref{Abel}). Since there is a one-to-one correspondence between pairs of polynomials $p,q$
and pairs of their primitives $P,Q$ defined above, we shall formulate all our results in terms of $P$ and $Q$. Below we always
assume that $Q$ with $Q(a)=Q(b)=0$ is fixed, while $P$ varies in the space ${\cal P}_d$ of all the polynomials of the degree up to 
$d$ vanishing at $a$ and $b$. 

Let us start with a description of the Composition Set $COS_{d,Q}$ of all the polynomials $P$ in ${\cal P}_d,$ such that $P$ and $Q$
satisfy the Composition Condition.

\bt \label{Int.cos} For $V \subset {\cal P}_d$ and for any polynomial $Q$ of degree at most $5$ the Composition
Set $COS_{d,Q}$ is a linear subspace in ${\cal P}_d$ of the dimension at most $[{d\over 2}]$. For $6 \leq \deg Q \leq 89$
the set $COS_{d,Q}$ is a union of at most two linear subspaces in ${\cal P}_d$, and for $\deg Q \geq 90$ the set $COS_{d,Q}$ is a 
union of at most three linear subspaces. The dimension of each of these subspaces is at most $[{d\over 2}]$, their double and triple
intersections have dimensions at most $[{d\over 6}]+1$ and $[{d\over 90}],$ respectively.
\et 
The rest of our results bound the dimension of the non-composition components, i.e. those which are not contained in $COS_{d,Q}$ (if
they exist).

\bt \label{Int1}
Consider equation (\ref{Abel}) with $Q$ fixed and $P$ varying in the space ${\cal P}_d$
of all the polynomials of the degree up to $d$ vanishing at $a$ and $b$. Then the dimension
of the non-composition components of the Center Set  of (\ref{Abel}), if they exist, does not exceed
$[{d\over 6}]+2.$ In particular, this dimension is of order at most one third of the maximal dimension of the 
Composition Center strata (which is of order $d\over 2$, being achieved on the compositions strata with the right factor 
$W(x)=(x-a)(x-b)$).
\et
The main steps in the proof of Theorem \ref{Int1} are the following: we consider the projective compactification $P{\cal P}_d$ of 
${\cal P}_d$, and use the fact, proved in \cite{broy}, that the Center Equations ``at infinity'' become the Moment equations. Therefore, 
to bound the dimensions of the affine non-composition components of the Center Set $CS$ in the {\it complex} affine space ${\cal P}_d$ 
it is enough to bound the dimensions of the non-composition components of the Moment vanishing set $MS$ ``at infinity'' in $P{\cal P}_d$. 
We show that these dimensions do not exceed $[{d\over 6}]+2,$ using a complete description of the Moments vanishing conditions, obtained 
in \cite{pak}.

More accurately, we define the set $ND$ of ``non-definite'' polynomials which provide non-composition solutions to the Moment equations,
and bound from above its dimension. Then the following theorem describes an inclusion structure at infinity of the sets we are interested
in: 

\bt \label{Int1.1}
For an algebraic set $Y\subset {\cal P}_d$ let $\bar Y$ denote the intersection of $Y$ with the infinite hyperplane of $P{\cal P}_d$. 
Then for each irreducible non-composition component $A$ of the affine Central Set $CS$ we have 
$\bar A \subset \bar {CS}\cap ND \subset \bar {MS} \cap ND$. Consequently, $\dim A \leq \dim (\bar {MS} \cap ND)+1$. 
\et

\smallskip

In many specific cases Theorem \ref{Int1.1} allows us to improve the general bound provided by Theorem \ref{Int1}.
In order to formulate corresponding results it is convenient to normalize points $a$ and
$b$ to be the points $-{{\sqrt 3}\over 2}$ and ${{\sqrt 3}\over 2}$, respectively.
Further, let ${\cal S} \subset {\cal P}$ be a subset of all polynomials $Q \in {\cal P}$
representable as a sum $Q=S_1(T_2)+S_2(T_3)$, where $S_1,S_2$  are arbitrary polynomials,
while $T_2,T_3$ are the Chebyshev polynomials of the degree $2$ and $3$, respectively
(notice that the normalization of the interval $[a,b]$ is chosen in such a way that $T_2(a)=T_2(b),$ $T_3(a)=T_3(b)$).
Below we show that the dimension of ${\cal S}\cap {\cal P}_d$ does not exceed $[{2\over 3}d]+1$,
so ``most'' of the polynomials $Q$ of degree $d$ cannot be represented in the above form.

\bt \label{Int2}
Let $P$ vary in the space ${\cal P}_9$. Then for each fixed
$Q \in {\cal P} \setminus {\cal S}$ the Center Set of (\ref{Abel}) consists of a Composition
Set and possibly a finite set of additional points. For an arbitrary fixed $Q$ the dimension
of the non-composition components of the Center Set of (\ref{Abel}) in ${\cal P}_9$ does not exceed one. For $P$
varying in the space ${\cal P}_{11}$ and for an arbitrary fixed $Q$ the dimension of the
non-composition components of the Center Set of (\ref{Abel}) does not exceed two.
\et
The next result heavily relies on computations with the second Melnikov coefficients.

\bt \label{Int3}
Let $P$ vary in the space ${\cal P}_9$. Then for each fixed $Q \in {\cal S}\cap {\cal P}_9$,
which is not a polynomial in $T_2$ or $T_3$, the Center Set
of (\ref{Abel}) consists of a Composition Set and possibly a finite set of additional points.
\et
Our last result (Theorem \ref{csinf8} in Section \ref{mainr} below) concerns the Center Set in subspaces of polynomials
with a special structure. Here we formulate its important particular case. Let $U_d$ consist of all polynomials
$P\in {\cal P}_d$ such that the degrees of $x,$ appearing in $P$ with the non-zero
coefficients are powers of prime numbers.
\bt \label{Int4}
Let $P$ vary in $U_d$.  Then for any fixed $Q$ the Center Set of (\ref{Abel}) in $U_d$
consists of a Composition Set and possibly a finite set of additional points.
\et

\section{Preliminaries: Poincar\'e mapping, Center Equations, and Composition
condition}
\setcounter{equation}{0}

\subsection{Poincar\'e mapping and Center Equations}\label{Cen.Eq}

Both the Center-Focus and the Smale-Pugh  problems can be naturally expressed in terms of the Poincar\'e
``first return'' mapping $y_b=G_{[a,b]} (y_a)$ along $[a,b]$. Let $y(x,y_a)$
denote the element around $a$ of the solution $y(x)$ of (\ref{Abel}) satisfying
$y(a)=y_a$. The Poincar\'e mapping $G_{[a,b]}$ associates to each initial value $y_a$ at $a$, sufficiently close to 
zero, the value $y_b$ at $b$ of the solution $y(x,y_a)$ analytically continued along $[a,b]$.

According to the definition above, the solution $y(x,y_a)$ is closed on $[a,b]$ if and only if $G_{[a,b]} (y_a)=y_a$. 
Therefore closed solutions correspond to the fixed points of $G_{[a,b]}$, and (\ref{Abel}) has a center if and only if
$G_{[a,b]}(y)\equiv y$. It is well known that $G_{[a,b]} (y)$ for small $y$ is given by a convergent power series

\be \label{eq1}
G_{[a,b]} (y)= y + \sum_{k=2}^\infty v_k(p,q,a,b)y^k.
\ee
Therefore the center condition $G_{[a,b]}(y)\equiv y$ is equivalent to an infinite sequence of algebraic equations on 
$p$ and $q$:

\be \label{eq2}
v_k(p,q,a,b)=0, \ k=2,3,\dots.
\ee
Each $v_k(p,q,a,b)$ can be expressed as a linear combination of certain iterated integrals of $p$ and $q$ along $[a,b]$ 
(see \cite{broy} and Theorem \ref{ceqinf} below).

\subsection{Projective setting and Center Equations at infinity over fixed $Q$}
Let ${\cal P}={\cal P}_{[a,b]}$ be the vector space of all complex polynomials $P$ satisfying
$P(a)=P(b)=0$, and ${\cal P}_d$ the subspace of $\cal P$ consisting of polynomials of degree at
most $d$. We always shall assume that the polynomials
\be \label{iiuu}
P(x)=\int_a^x p(\tau)d\tau, \ \ \ Q(x)=\int_a^x q(\tau)d,
\ee
defined above are elements of $\cal P$. This restriction is natural in the study of the center conditions since it is 
forced by the first two of the Center Equations \eqref{eq2}. Since \eqref{iiuu} provides a one to one correspondence between 
$(p,q)$ and $(P,Q)$, which is an isomorphism of the corresponding vector spaces, in order to avoid a cumbersome notation all 
the results below are formulated in terms of $(P,Q)$.

We shall assume that the points $a\ne b$ are fixed, and usually shall omit $a,b$ from the notations.

\medskip

From now on we shall assume that $Q \in {\cal P}_{d_1}$ is fixed, while $P$ varies in a certain linear subspace $V$ of the space 
${\cal P}_d$. This restrictive setting significantly simplifies the presentation, although it describes only ``slices'' of the 
Center Set. The approach of \cite{broy} and of the present paper can be extended to the full coefficients space of 
$(P,Q) \in {\cal P}_{d}\times {\cal P}_{d_1}$. We consider this extension as an important research direction, but it significantly
increases the complexity of the Algebraic Geometry involved, and is beyond the scope of the present paper. See \cite{broy} for a 
comparison of different possible settings of the problem.

\medskip

Let a subspace $V \subset {\cal P}$ be given. We shall consider the projective space $PV$ and the infinite hyperplane 
$HV \subset PV$. To construct $PV$ we introduce an auxiliary variable $\nu \in {\mathbb C}$ and consider the couples 
$(S,\nu), \ S\in V$, with $(S,\nu)$ and $(\lambda S,\lambda \nu)$ identified for any $\lambda \in {\mathbb C}, \ \lambda \ne 0$. 
The infinite hyperplane $HV$ is defined in $PV$ by the equation $\nu=0$.

Let us denote by ${\hat v}_k(p,q)={\hat v}_k(p,q,a,b)$ the ``homogenization" of the Center Equations $v_k(P,Q,a,b)=0$ with respect 
to the variable $P$. In other words, we multiply each term in $v_k$ by an appropriate degree of an auxiliary variable $\nu$ to make 
$v_k$ homogeneous.

\smallskip

Notice that the Center Equations can be considered in two ways: as polynomial equations in the coefficients of $P,Q$, or as symbolic
equations, containing ``symbolic iterated integrals'' of the form $\int  p \int q \int q ...$ (which can be interpreted as the 
poly-linear forms, i.e. polynomials, in the symbols $p,q$). Since each $p,q$ is a linear form in its coefficients, the degrees of the
polynomials in both interpretations are the same. Accordingly, the projective space $PV$ and the homogeneous ${\hat v}_k(p,q)=0$ can
be treated symbolically, till the moment where we have to actually integrate and get the explicit answer.  

\smallskip

We call ``Center Equations at infinity'' the restrictions of the homogeneous
Center Equations to the infinite hyperplane $HV$. They are obtained by
putting $\nu=0$ in the homogeneous equations described above. The following Theorem \ref{ceqinf} provides a description
of the Center Equations at infinity obtained in \cite{broy}. Take into account a different order of the polynomials
$p$ and $q$ in the Abel equation (\ref{Abel}) in the present paper and in \cite{broy}.

\bt \label{ceqinf} (\cite{broy})
For $k=2,4,\dots$ even, and $l={k\over 2}-1$ the Center Equations at infinity over $Q$
are given by vanishing of the generalized moments
\be \label{mominf}
v^\infty_k(P,Q)=m_l(P,Q)=\int_a^b P^l(x)q(x)dx=0.
\ee
For $k$ odd the Center Equations at infinity over $Q$ are given by vanishing of the coefficients of the
``second Melnikov function''
\be \label{melinf}
v^\infty_k(P,Q)=D_k(P,Q)=0,
\ee
represented by integer linear combinations $\sum n_\a I_\a$, with the sum running over all the iterated integrals in $p,q$
with exactly two appearances of $q.$ Here $\a = (\a_1,\ldots,\a_s),$ with exactly two of $\a_j$ equal to $1$, and the rest
equal to $2$, and with $\sum_{j=1}^s \a_j=k-1$. The integrals $I_\a$ are defined as

$$
I_\a=\int_a^b h_{\a_1}(x_1)dx_1 (\int_a^{x_1} h_{\a_2}(x_2)dx_2...(\int_a^{x_{s-1}} h_{\a_3}(x_s)dx_s)...),
$$
with $h_1=q, \ h_2=p$. The integer coefficients $n_\a$ are given as the products $n_\a=(-1)^s \prod_{r=1}^s (k-\sum_{j=1}^r \a_j)$.  
\et
In Proposition \ref{meln1} below the first four Melnikov equations at infinity $D_k(P,Q)=0$ are given
explicitly.

\subsection{Center, Moment, and Composition Sets}

Let us assume that $Q \in {\cal P}_{d_1}$ and a subspace $V \subset {\cal P}_d$
are fixed. We define the Center Set $CS=CS_{V,Q}$ as the set of $P\in V$
for which equation (\ref{Abel}) has a center. Equivalently, $CS$ is the set of $P\in V$
satisfying Center Equations (\ref{eq2}). The moment set $MS=MS_{V,Q}$ consists of
$P\in V$ satisfying Moment equations (\ref{mominf}).

To introduce Composition Set $COS=COS_{V,Q}$ we recall the polynomial
Composition Condition defined in \cite{bfy}, which is a special case of the general
Composition Condition introduced in \cite{al} (for brevity below we shall use
the abbreviation ``CC'' for the ``Composition Condition'').

\bd
Polynomials $P,Q$ are said to satisfy the ``Composition Condition''  on $[a,b]$
if there exist polynomials $\tilde P$, $\tilde Q$ and $W$ with $W(a)=W(b)$ such that
$P$ and $Q$ are representable as $$P(x)=\tilde P(W(x)), \ \  \ Q(x)=\tilde Q(W(x)).$$ The
Composition Set $COS_{V,Q}$ consists of all $P\in V$ for which $P$ and $Q$
satisfy the Composition Condition.
\ed
It is easy to see that the Composition Condition implies center for (\ref{Abel}), as well
as the vanishing of each of the moments and iterated integrals above. So we have
$COS \subset CS, \ COS \subset MS$.

\smallskip

Define $\bar {CS}, \bar {MS}, \bar {COS}$ as the intersections of the corresponding
affine sets with the infinite hyperplane $HV$. It follows directly from Theorem \ref{ceqinf} that
the following statement is true.

\bp \label{inclusion}
$\bar {COS} \subset \bar {CS} \subset \bar {MS}.$
\ep

Notice that $COS$ and $MS$
are homogeneous, and hence these sets are cones over $\bar {MS}, \bar {COS}$.
However, $CS$ a priori may be not homogeneous, and the connection of the
affine part $CS$ to $\bar {CS}$ may be more complicated.

Our main goal will be to compare the affine Center Set $CS$ with the Composition Set $COS$.
For this purpose we shall bound the dimension of the affine non-composition components
of $CS$, analyzing their possible behavior at infinity  (Sections 5,6). To obtain these
bounds we first describe the geometry of the Composition Set $COS$ (Section 3)
and compare the Moment set $MS$ and its subset $COS$ (Section 4).

\section{The structure of the Composition Set}
\setcounter{equation}{0}

Geometry of the Composition Set reflects the algebraic structure of
polynomial compositions, which is well known to provide rather subtle
phenomena. In comparison with the classical theory developed by Ritt (\cite{rit})
we are interested in what we call below $[a,b]$-compositions, i.e. compositions of
polynomials under requirement that some of the factors take equal values at the
points $a$ and $b$.

\subsection{Elements of Ritt's theory}

Let us recall first some basic facts on polynomial composition algebra, including
the classical first and second Ritt theorems (\cite{rit}).

\bd A polynomial $P$ is called indecomposable if it cannot be represented
as $P(x)=R\circ S(x)=R(S(x))$ for polynomials $R$ and $S$ of degree greater than
one. A decomposition $P=P_1\circ P_2\circ \dots \circ P_r$ is called maximal
if all $P_1,\dots,P_r$ are indecomposable and of degree greater than one. Two decompositions
$P=P_1\circ P_2\circ \dots \circ P_r$ and $P=Q_1\circ Q_2\circ \dots \circ Q_r,$
maximal or not, are called equivalent (notation ``$\sim$'') if there exist polynomials of
degree one $\mu_i, \ i=1,\dots,r-1,$ such that
$P_1= Q_1\circ \mu_1, \ P_i=\mu_{i-1}^{-1}\circ Q_i\circ \mu_i, \ i=2,\dots,r-1,$ and
$P_r=\mu_{r-1}^{-1}\circ Q_r$.
\ed

The first Ritt theorem (\cite{rit}) states that any two maximal decompositions of a polynomial
$P$ have an equal number of terms, and can be obtained from one another by a sequence
of transformations replacing two successive terms $A\circ C$ with $B\circ D$, such
that

\be \label{riii}
A\circ C = B\circ D.
\ee
Let us mention that decompositions of a polynomial $P$ into a composition of two polynomials,
up to equivalence, correspond in a one-to-one way to imprimitivity systems of the monodoromy
group $G_P$ of $P$ (see e.g. \cite{rit} or \cite{pak3}). In their turn imprimitivity systems of $G_P$
are in a one-to-one correspondence with subgroups $A$ of $G_P$ containing the stabilizer $G_\omega$
of a point $\omega \in G$.
In particular, for a given
polynomial $P$ the number of its right composition factors $W$, up to the change
$W\rightarrow \lambda\circ W$, where $\lambda$ is a polynomial of degree one, is
finite. Below we shall call (with a slight abuse of notation) two right composition
factors $W$ and $\l\circ W$ of $P$, where $\lambda$ is a polynomial of degree one,
equivalent, and write $W \sim \l\circ W$. We also usually shall write just ``right
factor''  of $P$ instead of ``compositional right factor''.

\smallskip

The first Ritt theorem reduces the description of maximal decompositions of
polynomials to the description of indecomposable polynomial solutions of the
equation (3.1). It is convenient to start with the following result (\cite{en}):
if polynomials $A,B,C,D$ satisfy \eqref{riii}, then there exist polynomials
$U,V,\hat A, \hat B, \hat C, \hat D,$ where
\be \label{nach} \deg U=\gcd(\deg A, \deg B), \ \ \ \deg V=\gcd(\deg C, \deg D),\ee such that
\be \label{xcv}
A=U\circ \hat A, \ \ \ B=U\circ \hat B, \ \ \ C= \hat C\circ V, \  \ \ D= \hat D\circ V,
\ee
and \be \label{riii+} \hat A \circ \hat C = \hat B \circ \hat D.\ee In particular,
if $\deg A = \deg B$, then necessarily $A\circ C$ and $B\circ D$ are equivalent as decompositions.
More generally, if $\deg B\vert \deg A$, then there exists a polynomial $W$ such that the
equalities $$A=B\circ W, \ \ \  D=W\circ C$$ are satisfied.

Note that the above result concerning the reduction of \eqref{riii} to \eqref{riii+}
is equivalent to the statement that the lattice of imrimitivity systems of the monodromy
group $G$ of a polynomial $P$ of degree $n$ is isomorphic to a sublattice of the lattice $L_n$
consisting of all divisors of $n$, where by definition $$d_1\wedge d_2=\GCD(d_1,d_2), \ \ \ d_1\vee d_2=\rm{LCM} (d_1,d_2)$$
(see \cite{mp2}). For example, for the
polynomials $z^n$ the corresponding lattices consist of all divisors of $n$ since for any $d\vert n$ the equality
$z^n=z^d\circ z^{n/d}$ holds. The same is true for the Chebyshev polynomials $T_n$ since the equality
$T_n(\cos \phi)=\cos n\phi$ implies that $T_n=T_d\circ T_{n/d}$ for any $d\vert n.$
On the other hand, for an indecomposable polynomial
$P$ the corresponding lattice contains only elements $1$ and $n.$

\medskip

The second Ritt theorem (\cite{rit}) states that if $A,B,C,D$ satisfy \eqref{riii} and degrees of $A$ and $B$ as well as of
$C$ and $D$ are coprime, then there exist linear polynomials $U,V$ such that \eqref{xcv} and \eqref{riii+} hold, and,
up to a possible replacement of
$\hat A$ by $\hat B$ and $\hat C$ by $\hat D$, either

\be \label{sol1}
\hat A \circ \hat C \sim z^n\circ z^r R(z^n), \ \ \ \ \hat B \circ \hat D \sim z^r R^n(z)\circ z^n,
\ee

\noindent where $R(z)$ is a polynomial, $r\geq 0, n \geq 1,$ and $\GCD (n,r)=1,$ or

\be \label{sol2}
\hat A \circ \hat C \sim T_n\circ T_m, \ \hat B \circ \hat D \sim T_m\circ T_n,
\ee

\noindent where $T_n$ and $T_m$ are the  Chebyshev polynomials, $n,m\geq 1$, $\GCD (n,m)=1$.
In particular, this holds when $A,B,C,D$ solving (3.1) are indecomposable, and the
decompositions $A\circ C$ and $B\circ D$ are non-equivalent, since in this case the degrees of
polynomials $U,V$ in \eqref{nach} and \eqref{xcv} are necessarily equal to one.

Clearly,  the second Ritt theorem together with the previous result imply the following statement:
if $A,B,C,D$ satisfy \eqref{riii}, then there exist  polynomials $U,V$ such that \eqref{nach},
\eqref{xcv}, \eqref{riii+} hold, and, up to a possible replacement of
$\hat A$ by $\hat B$ and $\hat C$ by $\hat D$, either \eqref{sol1} or \eqref{sol2} holds.

\subsection{$[a,b]$-Compositions}

Now we return to $[a,b]$-compositions, i.e. compositions of polynomials under the
requirement that some of the right factors take equal values at two distinct points
$a$ and $b$.

\bd \label{abcomp} Let a polynomial $P$ satisfying $P(a)=P(b)$ be given. We call
polynomial $W$ a right $[a,b]$-factor of $P$ if $P={\tilde P}\circ W$ for some polynomial $\tilde P,$
and $W(a)= W(b)$. A polynomial $P$ is called $[a,b]$-indecomposable, if $P(a)=P(b)$ and $P$ does not
have right $[a,b]$-factors non-equivalent to $P$ itself.
\ed

\noindent{\bf Remark.} Notice that any right $[a,b]$-factor of $P$ necessary has degree greater than
one, and that $[a,b]$-indecomposable $P$ may be decomposable in the usual sense.

\bp \label{iq} Any polynomial $P$ up to equivalence has a finite number of $[a,b]$-indecomposable
right factors $W_j, \ \ j=1,\dots, s$. Furthermore, each right $[a,b]$-factor $W$ of $P$ can be
represented as $W=\tilde W(W_j)$ for some polynomial $\tilde W$ and $j=1,\dots, s.$

\ep \pr As it was mentioned above up to equivalence
there are only finitely many general right factors $W$ of $P$. In particular,
this is true for $[a,b]$-indecomposable right $[a,b]$-factors $W_j$ of $P$.

Now let $W$ be a right $[a,b]$-factor of $P$. If it is $[a,b]$-indecomposable,
then by the first part of the proposition $W=\l\circ W_j$ for some $j=1,\dots, s$, with $\l$ a linear
polynomial. Otherwise, $W$ can be represented as $W=V\circ \hat W$, where
$\hat W$ is a right $[a,b]$-factor of $P$ and $\deg V > 1$. Since
$\deg \hat W < \deg W$, it is clear that continuing this process we ultimately will find an
$[a,b]$-indecomposable right factor $W_j$ of $P$ such that $W=\tilde W(W_j)$. $\Box$

\medskip

An easy consequence of Proposition \ref{iq} is the following description of the Composition
Set given in \cite{broy}:

\bp \label{bry} Let $W_j,$ $j=1,\dots, s,$ be all indecomposable right $[a,b]$-factors
of $Q$. Then the set $COS_{V,Q}$ is a union of the linear subspaces
$L_j\subset V$, $j=1,\dots, s$, where $L_j$ consists of all the polynomials
$P \in V$ representable as $P=\tilde P(W_j), \ \ j=1,\dots, s,$ for a certain
polynomial $\tilde P$.
\ep
It has been recently shown in \cite{pak} that for any $P\in {\cal P}$ the number
$s$ of its non-equivalent $[a,b]$-indecomposable right factors can be at
most three. Moreover, if $s>1$ then these factors  necessarily
have a very special form, similar to what appears in Ritt's description above.

The precise statement is given  by the following theorem (\cite{pak}, Theorem 5.3):

\bt \label{three} Let complex numbers $a\ne b$ be given. Then for any polynomial
$P\in {\cal P}_{[a,b]}$ the number $s$ of its $[a,b]$-indecomposable right
factors $W_j$, up to equivalence, does not exceed 3.

Furthermore, if $s=2$, then either

$$P=U\circ z^{rn}R^n(z^n)\circ U_1, \ W_1=z^n\circ U_1,
\ W_2=z^r R(z^n)\circ U_1,$$ where
$R, U, U_1$ are polynomials, $r>0, n>1,$ $\GCD (n,r)=1,$
or

$$ P=U\circ T_{nm}\circ U_1, \ W_1=T_n\circ U_1, \ W_2=T_m\circ U_1,$$
where
$U, U_1$ are polynomials,
$n,m>1, \ \gcd (n,m)=1.$

On the other hand, if
$s=3$ then
$$P=U\circ z^2 R^2(z^2)\circ T_{m_1m_2}\circ U_1,$$ $$W_1=T_{2m_1}\circ U_1, \ \
\ W_2=T_{2m_2}\circ U_1, \ \ \  W_3=zR(z^2)\circ T_{m_1m_2}\circ U_1,$$
where
$R, U, U_1$ are polynomials,
$m_1,m_2 >1$
are odd, and $\gcd (m_1,m_2)=1.$

\et
Notice that
in all the cases above $U_1(a)\ne U_1(b)$ while $W_j(a)=W_j(b)$.

\smallskip
We are interested in the stratification of the space ${\cal P}_d$ of polynomials $P$ of
degree $d$ according to the structure of their $[a,b]$-indecomposable right $[a,b]$-factors.
Following Theorem \ref{three} let us use the following notations for the appropriate strata:

\bd \label{dec} Let $DEC^d_s(a,b)\subset {\cal P}_d$ denote the set of polynomials $P$ of degree
at most $d$ satisfying $P(a)=P(b)=0$ and possessing exactly $s$
non-equivalent $[a,b]$-indecomposable right factors. For $s=1$ we write
$DEC^d_1(a,b)= DEC^d_{1,0}(a,b) \cup DEC^d_{1,1}(a,b).$ Here $DEC^d_{1,0}(a,b)$ consists of polynomials $P$ for which 
their only indecomposable right factor $W$ is equivalent to $P$. In turn, $DEC^d_{1,1}(a,b)$ consists of $P$ for which 
$W$ is not equivalent to $P$, and hence $\deg W < \deg P$.
\ed

As a first consequence of Theorem \ref{three} we get upper bounds on the dimensions of the sets
$DEC^d_s(a,b)$ considered as subsets of the complex space $\mathbb C^{d-1}$, which we identify with
${\cal P}_d$.

\bp \label{propp} $DEC^d_{1,0}(a,b)$ consists of $[a,b]$-indecomposable polynomials $P\in {\cal P}_d$,
and its dimension is $d-1$. We have $DEC^d_{1,1}(a,b)=\emptyset$ for $d\leq 3$, and
$\dim DEC^d_{1,1}(a,b) \leq [{d\over 2}]$ for $d\geq 4$. $DEC^d_2(a,b)=\emptyset$ for $d\leq 5$, and
$\dim DEC^d_2(a,b) \leq [{d\over 6}]+1$ for $d\geq 6$. $DEC^d_3(a,b)=\emptyset$ for $d\leq 89$, and
$\dim DEC^d_3(a,b) \leq [{d\over 90}]$ for $d\geq 90$.
\ep
\pr
Assume we are given $l$ parametric families of polynomials \ \ \
${\cal S}_r=\{S_r(\tau_r,z)\}, \ r=1,...,l$, with $\tau_r \in \Sigma_r \subset {\mathbb C}^{n_r}$
being the parameters of ${\cal S}_r$. We assume that the degree of the polynomials $S_r(\tau_r,z)$
remains constant and equal to $d_r$ for all the values of the parameters $\tau_r \in \Sigma_r$. Put
$\tau=(\tau_1,\dots,\tau_l),$ and let
$$
P_\tau=S_1(\tau_1) \circ S_2(\tau_2) \circ \dots \circ S_l(\tau_l).
$$
The degree of the polynomials $P_\tau$ of this form is $d_1\cdot ... \cdot d_l$ and they form
a parametric family with the parameters $\tau=(\tau_1, ... \tau_l)\in {\mathbb C}^n,$ where
$n=n_1+...+n_l$.

\smallskip

The dimension $D$ of the stratum $S$ in $\cal P$ formed by the polynomials $P_\tau$ as above is at
most $n$, and it may be strictly less than $n$ since the parametric representation as above may
be redundant. The requirement $P_\tau \in {\cal P}_d$ is equivalent to $d_1\cdot ... \cdot d_l \leq d.$

So in order to bound from above the dimensions of the strata $DEC^d_s(a,b),$ we have to accurately estimate
the number $D \leq n_1+...+n_l$ of free parameters, and the degrees $d_1, ... , d_l$
in composition representations of the corresponding polynomials $P$, provided by Theorem \ref{three}. We have to
take into account the redundancy in the parametric representation, and then to maximize $D$ under the constraint 
$d_1\cdot ... \cdot d_l \leq d.$

\smallskip

Notice that ${\cal P}_d=\cup_{s=1}^3 DEC^d_s(a,b)$. Let us now consider the sets $DEC^d_s(a,b)$ for $s=1,2,3$ case 
by case. We shall see below that all the strata $DEC^d_s(a,b),$ besides the stratum $DEC^d_{1,0}(a,b),$ consisting 
of $[a,b]$-indecomposable polynomials $P$, have dimension strictly smaller than $\dim {\cal P}_d = d-1.$ Hence
$\dim DEC^d_{1,0}(a,b)=d-1.$ (This follows immediately also from the fact that $DEC^d_{1,0}(a,b)$ consists of 
{\it generic} polynomials in ${\cal P}_d$). 

Now, each $P \in DEC^d_{1,1}(a,b)$ has a form $P=S_1\circ S_2$,
with $\deg S_1=d_1 >1, \deg S_2=d_2>1,$ since we assume that $P$ possesses a right
$[a,b]$-factor $S_2$, not equivalent to $P$. In this case $d\geq d_1d_2$ is at least $4$, and
$S_1$ and $S_2$ can be any polynomials of degrees $d_1$ and $d_2$ with the only restrictions
$S_2(a)=S_2(b)$ and $S_1(S_2(a))=0$. Hence $n_1=d_1, \ n_2=d_2$. On the space ${\mathbb C}^{n_1+n_2}$
of the parameters of $(S_1,S_2)$ acts a two-dimensional group $\Gamma$ of linear polynomials $\gamma$.
It acts by transforming $(S_1,S_2)$ into $(S_1\circ \gamma, \gamma^{-1}S_2).$ This action preserves $P$.
Accordingly, we have to maximize $D=d_1+d_2-2$ under the constrain $d_1d_2\leq d$. For $d$ even this
maximum is achieved for $d_1=2$ or $d_2=2$ and it is ${{d}\over 2}$. For $d$ odd still $d_1=2$ or
$d_2=2$, but the maximum of $D$ is ${{d-1}\over 2}$. Finally we get
$\dim DEC^d_1(a,b) \leq [{d\over 2}].$

\smallskip

Now let us consider the case $s=2$. In this case by Theorem \ref{three} we have two options.

\smallskip

The first option is that $P=U\circ z^{rn}R^n(z^n)\circ U_1, $ where $U(z), R(z), U_1(z)$
are polynomials, $r>0, n>1,$ and $\gcd (n,r)=1,$ and $z^n$ and $z^r R(z^n)$ take equal
values at $U_1(a) \ne U_1(b)$.

\smallskip

Here, denoting the degrees of $U,U_1,R$ by $k,m,l\geq 1,$ respectively, we get
$\deg P= k\cdot n(r+ln)\cdot m \geq 6,$ while the number of the independent parameters,
i.e. the dimension of the corresponding strata is at most $k+l+m-1$ (we take into account the
requirements $W_1(a)=W_1(b), \ W_2(a)=W_2(b), \ P(a)=P(b)=0,$ and the fact that the scaling
parameters of $U$ and of $R$ act equivalently on $P$). So we have to maximize $k+l+m-1$ under
the constrain $k\cdot n(r+ln)\cdot m \leq d.$ The variables are integers
$k\geq 1,l\geq 1,m\geq 1, r\geq 1, n\geq 2, \ \gcd (n,r)=1$.

\smallskip

Let us first fix $l,r,n$. As above, the maximum of $k+l+m-1$ is
attained either for $k=1, m=[{d\over {n(r+ln)}}],$ or for $k= [{d\over {n(r+ln)}}], m=1.$
In both cases it is $l+ [{d\over {n(r+ln)}}],$ and this expression increases as $l$
decreases. So we can put $l=1$ and so we get $[{d\over {n(r+n)}}]+1.$ Once more, this
expression increases as $n,r$ (which do not enter the maximized sum) decrease. Their minimal
possible values are $r=1,n=2$ and we get $k+l+m-1=[{d\over 6}]+1.$
\smallskip

The second option is that $P=U\circ T_{nm}\circ U_1,$ with $n,m>1, \ \gcd (n,m)=1,$ and
$T_m$ and $T_n$ take equal values at $U_1(a)$ and $U_1(b)$. Denote the degrees of $U$ and
$U_1$ by $k$ and $l$, respectively. We get $\deg P=klmn\geq 6$, while the number of the
independent parameters, i.e. the dimension of the corresponding strata, is at most $k+l-1$
(we take into account the requirements that $T_m$ and $T_n$ take equal values at $U_1(a)$
and $U_1(b)$, and $P(a)=P(b)=0$). By exactly the same reasoning as above we conclude that the maximal
dimension of the corresponding strata is achieved as either $\deg U=1$ or $\deg U_1=1$, and
it is at most $[{d\over {mn}}].$ The minimal possible values for $m,n$ here are $2$ and $3$,
so we get the bound $[{d\over 6}]$ which is smaller than the one above.

\smallskip
\smallskip

It remains to consider the case $s=3$. In this case by Theorem \ref{three} we have
$P=U\circ z^2 R^2(z^2)\circ T_{m_1m_2}\circ U_1,$ with $U,R,U_1$ as above, $m_1,m_2 >1$
odd, and $\gcd (m_1,m_2)=1.$ In addition, $T_{2m_1}, T_{2m_2}$ and $zR(z^2)\circ T_{m_1m_2}$
take equal values at $U_1(a) \ne U_1(b)$.

As above, denoting the degrees of $U,U_1,R$ by $k,m,l,$ respectively, we get
$\deg P = k\cdot (4l+2)m_1m_2\cdot m \geq 90.$ The number of the independent parameters,
i.e. the dimension of the corresponding strata, is here at most $k+l+m-2$ (we take into account,
besides the requirements that $W_1,W_2,W_3$ take equal values at $a,b$, and $P(a)=P(b)=0$, also
the fact that the scaling parameters of $U$ and of $R$ act equivalently on $P$). Maximizing the
last expression exactly as above, we conclude that the maximum is achieved for
$l=1, m_1=3, m_2=5,$ and either $k=1, m=[{d\over {(4l+2)m_1m_2}}]=[{d\over 90}],$ or
$m=1, k= [{d\over 90}].$ This maximum is equal to $[{d\over {90}}].$ This completes the proof
of Proposition \ref{propp}. $\square$.

\medskip

Based on Proposition \ref{propp}  and Theorem \ref{three} we can now give a much more accurate
description of the Composition Set $COS_{V,Q}$ for $V \subset {\cal P}_d$:

\bt \label{cos} For $V \subset {\cal P}_d$ and for any polynomial $Q$ of degree at most $5$ the Composition
Set $COS_{V,Q}$ is a linear subspace $L$ in $V$ with $\dim L \leq [{d\over 2}]$. For $6 \leq \deg Q \leq 89$
the set $COS_{V,Q}$ is a union of at most two linear subspaces in $V$, and for
$\deg Q \geq 90$ the set $COS_{V,Q}$ is a union of at most three linear subspaces. The
dimension of each of these subspaces is at most $[{d\over 2}]$, their double and triple
intersections have dimensions at most $[{d\over 6}]+1$ and $[{d\over 90}],$ respectively.
\et
\pr It is sufficient to consider the case $V = {\cal P}_d$. Let $W_j, \ j=1,\dots,s$ be all the mutually
prime right $[a,b]$-factors of $Q$. By
Proposition \ref{propp}, for $Q$ of degree at most $5$ we have $s=1$. For $6 \leq \deg Q \leq 89$
we have $s\leq 2$ and for $\deg Q \geq 90$ we have $s\leq 3.$ Next, by Proposition \ref{bry},
$COS_{{\cal P}_d,Q}$ is a union of linear subspaces $L_j = \{P\in {\cal P}_d, \ P=\tilde P(W_j)\}$.

Next notice that if $\deg W_j = d,$ then $L_j$ is one-dimensional, and if $\deg W_j < d,$ then
$L_j \subset DEC^d_{1,1}(a,b)\cup DEC^d_2(a,b) \cup DEC^d_3(a,b).$ We also have
$L_i\cap L_j \subset DEC^d_2(a,b) \cup DEC^d_3(a,b), \ L_i\cap L_j \cap L_k \subset DEC^d_3(a,b).$ All the required
bounds on the dimensions of $L_j$ now follow directly from Proposition \ref{propp}. $\square$

\medskip

\noindent {\bf Remark}. In fact, the dimensions of the linear subspaces $L_j$ and of their
intersections may be strongly smaller than the bounds in Theorem \ref{cos}. The reason is that
in this theorem we do not take into account, for example, the fact, that if $Q$ has mutually
prime right $[a,b]$-factors $W_1,W_2$, then their degrees, by Theorem \ref{three}, cannot both be
equal to two. Another reason is that in the setting of Theorem \ref{cos} the right factors are
fixed, while in Proposition \ref{propp} they are variable, which also decreases the dimensions
of the strata of $COS_{{\cal P}_d,Q}$ in comparison with the strata $DEC^d_s(a,b).$

\section{Moment vanishing versus Composition}\label{Mom.Comp}
\setcounter{equation}{0}

The main result of \cite{pak} can be formulated as follows:

\bt \label{yui} Let $P$ with $P(a)=P(b)$ be given, and let $W_j,$ $j=1,\dots,s,$
be all its non-equivalent $[a,b]$-indecomposable right $[a,b]$-factors. Then
for any polynomial $Q$ all the moments $m_k=\int_a^b P^k(x)q(x)dx,$ $k\geq 0,$
vanish if and only if
$Q=\sum_{j=1}^s Q_j,$ where $Q_j=\tilde Q_j(W_j)$ for some polynomial $\tilde Q_j.$
\et
This theorem combined with Theorem \ref{three} provides an explicit description for
vanishing of the polynomial moments. In order to use it for the study of the Moment
Set, let us introduce the notions of ``definite'' and ``co-definite'' polynomials.

\bd \label{def.codef} Let $V,V_1 \subset {\cal P}={\cal P}_{[a,b]}$ be fixed linear
spaces. A polynomial $P \in {\cal P}$ is called $V_1$-definite if for any polynomial
$Q\in V_1$ vanishing of the moments
$m_k=\int_a^b P^k(x)q(x)dx,$ $k\geq 0,$ implies Composition Condition on $[a,b]$ for
$P$ and $Q$. The set of such $P$ is denoted $D_{V_1}$.

A polynomial $Q \in {\cal P}$ is called $V$-co-definite if for any polynomial $P\in V$
vanishing of the moments $m_k=\int_a^b P^k(x)q(x)dx,$ $k\geq 0,$ implies Composition
Condition on $[a,b]$ for $P$ and $Q$. The set of such $Q$ is denoted $COD_V$.

If $V_1 ={\cal P}$ or $V= {\cal P}$ (with respect to the corresponding $P$ or $Q$) 
we call polynomials defined above $[a,b]$-definite
or $[a,b]$-co-definite correspondingly, and denote their sets by $D$ or $COD$.
\ed
Definite polynomials have been initially introduced and studied in \cite{pry}.
Some their properties have been described in \cite{pp}. The notion of a co-definite
polynomials is apparently new (although some examples have appeared in \cite{broy}).
Below we give a characterization of definite and co-definite polynomials, but many
questions still remain open.

\subsection{Definite polynomials}

Theorem \ref{yui} allows us to give a complete description of $[a,b]$-definite
polynomials:

\bt \label{definite} A polynomial $P$ is $[a,b]$-definite if and only if it has,
up to equivalence, exactly one $[a,b]$-indecomposable right factor $W$.
\et
\pr
Assume that $P$ has exactly one $[a,b]$-indecomposable right factor $W$. By Theorem
\ref{yui} for any polynomial $Q$ vanishing of $m_k$ for all $k\geq 0$ implies that
there exist $\tilde Q$ such that $Q=\tilde Q(W)$, so Composition Condition on $[a,b]$
is satisfied for $P$ and $Q$. Hence, by Definition \ref{def.codef}, $P$ is $[a,b]$-definite.

Assume now that  $P$ has two non-equivalent $[a,b]$-indecomposable right factors $W_1$, $W_2$,
and show that the solution $Q=W_1+W_2$ cannot be represented in the form $Q=\tilde Q(W)$,
where $W$ is an $[a,b]$-right factor of $P$ and $\tilde Q$ is a polynomial (cf. \cite{acoun}).
First observe that $W_1$ and $W_2$ have different degrees for otherwise equalities
\eqref{xcv} imply that $W_1$ and $W_2$ are equivalent. Thus, without loss of generality
we may assume that $\deg W_2>\deg W_1$, and so $\deg Q=\deg W_2,$ implying that if
$Q=\tilde Q(W)$ then $\deg W\vert \deg W_2$. Therefore, using \eqref{xcv} again, we conclude
that $W_2=U(W)$ for some polynomial $U$. Furthermore, if $\deg W<\deg W_2$, then we obtain a
contradiction with the assumption that $W_2$ is an $[a,b]$-indecomposable right factor of $P.$
On the other hand, if $\deg W=\deg W_2$, then as above we conclude that $W$ and $W_2$ are
linear equivalent implying that $W_1=Q-W_2$ is a polynomial in $W_2$ in contradiction
with the assumption  $\deg W_2>\deg W_1$. $\square$

\medskip

Corollaries \ref{cc1}-\ref{cc2} below  were proved  in \cite{pp}. Here we give another
proof of these results basing on Theorem \ref{definite} and the second Ritt theorem. We
believe that these ``more algebraic'' proofs clarify to some extent the structure of
definite polynomials, which still presents a lot of open questions (see \cite{pry}).
We also extend a classification of non-definite polynomials whose degree does not exceed
nine, given in \cite{pp}, up to degree eleven.

\bc \label{cc1} Let $p$ be a prime. Then each polynomial $P$ of degree  $p^s$, $s\geq 1,$ is
$[a,b]$-definite for any $a,b\in \C, \ a\ne b.$
\ec
\pr Indeed, since imprimitivity systems of $G_P$ form a sublattice of $L_{p^s}$ (see definition
on page 9), if $W_1$, $W_2$ are arbitrary right factors of $P$, then either
$W_1$ is a polynomial in $W_2$ or $W_2$ is a polynomial in $W_1$. Therefore, such $P$ can
not have two non equivalent $[a,b]$-indecomposable right factors. $\square$

\bc \label{cc2}
If at least one of points $a$ and $b$ is not a critical point of a polynomial $P$,
then $P$ is $[a,b]$-definite.
\ec
Assume that $P$ is not $[a,b]$-definite and let $W_1$, $W_2$ be its non linear equivalent
$[a,b]$-indecomposable right factors. Then the second Ritt theorem implies that there exist
polynomials of degree one $\mu_1$, $\mu_2$ and polynomials $U$, $W$ such that either
\be \label{first} P=U\circ z^{rs}R^n(z^n)\circ W, \ \ \ W_1=\mu_1\circ z^n \circ W,
\ \ \ W_2=\mu_2\circ z^sR(z^n) \circ W,
\ee
where $R$ is a polynomial and $GCD(s,n)=1$, or
\be \label{second} P=U\circ T_{nm}\circ W, \ \ \ W_1=\mu_1\circ T_n \circ W,
\ \ \ W_2=\mu_2\circ T_m \circ W,
\ee
where $T_n,$ $T_m$ are the Chebyshev polynomials and $GCD(n,m)=1.$ Furthermore, since
$W_1$, $W_2$ are $[a,b]$-indecomposable and non equivalent, the inequality $W(a)\neq W(b)$
holds. In particular, $n>1,$ since $W_1(a)=W_1(b)$.

It is easy to see that if \eqref{first} holds, then the equalities
$$
W_1(\tilde a)=W_1(\tilde b), \ \ \  W_2(\tilde a)=W_2(\tilde b),
$$ where
$$
\tilde W_1=z^n, \ \ \  \tilde W_2=z^sR(z^n), \ \ \ \tilde a=W(a), \ \ \ \tilde b=W(b),
$$
taking into account the equality $GCD(s,n)=1,$ imply that the number $\tilde a^n=\tilde b^n$
is a root of the polynomial $R$. It follows now from the first formula in \eqref{first} by
the chain rule that both $a$ and $b$ are critical points of $P.$

If \eqref{second} holds, then, taking into account the identity \be \label{fifth} T_l\circ
\frac{1}{2}\left(z+\frac{1}{z}\right)=\frac{1}{2}\left(z+\frac{1}{z}\right)\circ z^l\ee and
the eqaulity $GCD(m,n)=1$, it is easy to see that there exist $\alpha,\beta\in \C$ such that
\be \label{fourth}\tilde  a=\frac{1}{2}\left(\alpha+\frac{1}{\alpha}\right), \ \ \
\tilde  b=\frac{1}{2}\left(\beta+\frac{1}{\beta}\right), \ \ \
\alpha^n=\beta^n,\  \ \ \alpha^m=\frac{1}{\beta^m},\ee
where as above $\tilde  a=W(a),$ $\tilde  b=W(b).$
Furthermore, $\alpha^2\neq 1$. Indeed, otherwise the equalities
$${\bar\alpha}^n={\bar \beta}^n,\  \ \ {\bar\alpha}^m=\frac{1}{{\bar\beta}^m},$$
where $\bar \alpha =\alpha^2,$ $\bar \beta =\beta^2$, taking into account
the equality $GCD(m,n)=1$, imply that $\beta^2=1$. Since $\tilde a\neq \tilde b$ this yields that either
$\tilde a=-1,$ $\tilde b=1,$ or $\tilde a=1,$ $\tilde b=-1.$ On the other hand, since $GCD(m,n)=1$,
without loss of generality we may assume that $m$ is odd implying that
$T_m(\tilde a)\neq T_m(\tilde b)$ for such $\tilde a$ and $\tilde b$ since $T_m(-1)=-1,$ $T_m(1)=1.$
Similarly,  $\beta^2\neq 1$. Finally, observe that equalities \eqref{fourth} yield that
$\alpha^{mn}=\pm 1,$ $\beta^{mn}=\pm 1,$
implying that \be \label{cvb} T_{mn}(\alpha)=\pm 1,\ \  \ T_{mn}(\beta)=\pm 1.\ee

In order to finish the proof observe that the equality $T_n(\cos \phi)=\cos n\phi$ implies easily that
the polynomial $T_n$ has exactly two critical values $\pm 1$ and that the only points in the preimage
$T_n^{-1}\{\pm 1\}$ which are not critical points of $T_n$ are the points $\pm 1$. Therefore,
the equalities \eqref{cvb} taking account that
$\alpha\neq \pm 1,$ $\beta\neq \pm 1$  imply that $\alpha$ and $\beta$ are critical points of $T_{mn}$
and hence critical points of $P$ by the chain rule. $\square$

\medskip

Theorem \ref{definite} combined with the second Ritt theorem allows us, at list in principle,
to describe explicitely all the non-definite polynomials up to a given degree. In particular,
the following statement holds:

\bt \label{cla} For given $a\ne b$ non-definite polynomials $P \in {\cal P}_{11}$
appear only in degrees $6$ and $10$ and have, up to change $P\rightarrow \lambda \circ P,$
where $\lambda$ is a polynomial of degree one, the following form:

\smallskip

1. $P_6=T_6 \circ \tau$, where $T_6$ is the Chebyshev polynomial of degree $6$ and $\tau$
is a polynomial of degree one transforming $a,b$ into
$-{{\sqrt 3} \over 2}, {{\sqrt 3} \over 2}$.

\smallskip

2. $P_{10}=z^2 R^2(z^2)\circ \tau$, where $R(z)= z^2 + \gamma z + \delta$ is an arbitrary
quadratic polynomial satisfying $R(1)=0$ i.e. $\gamma+\delta=-1$, and $\tau$ is a
polynomial of degree one transforming $a,b$ into $-1,1$.

\et
\pr First of all observe that if in Ritt's second theorem (Section 3.1 above) the degree of one
of polynomials satisfying \eqref{riii+} is two, then solutions \eqref{sol2} may be written
in form \eqref{sol1}. Indeed, for odd $n$ the equality
\be \label{chh}
T_n(z)=zE_n(z^2)
\ee
holds for some polynomial $E_n$. Furthermore, $T_2=\theta\circ z^2,$ where $\theta=2z-1$,
and hence
$$
zE_n(z^2)\circ \theta\circ z^2=T_n\circ T_2=T_2\circ T_n=\theta\circ T_n^2=\theta
\circ zE^2_n(z)\circ z^2.
$$
Since the last equality implies the equality
$$
zE_n(z^2)\circ \theta=\theta \circ zE^2_n(z),
$$
we conclude that
\be \label{tog}
T_n =\theta \circ zE^2_n(z)\circ \theta^{-1}, \ \ \ T_{2n} =\theta \circ z^2E^2_n(z^2).
\ee
Therefore, the equality
$$
T_n\circ T_2=T_2\circ T_n
$$
may be written in the form
\be \label{xomiak2}
(\theta \circ zE^2_n(z)\circ \theta^{-1})\circ (\theta\circ z^2)=(\theta\circ z^2)\circ
zE_n(z^2).
\ee
Now we are ready to prove the theorem.

Since each integer $i$, $2\leq i <  11$, distinct from 6 or 10 is either a prime or a power of
a prime, it follows from Corollary \ref{cc1} that $P$ is $[a,b]$-definite unless $\deg P=6$ or
$\deg P=10$. It follows now from the second Ritt theorem and the remark above that if $\deg P=10$,
then $P$ has the form given above. Similarly, if $\deg P=6$, then $P=z^2 R^2(z^2)\circ \tau$, where
$R$ is a polynomial satisfying $R(1)=0$. However, since in this case the degree of $R$ equals one,
up to change $P\rightarrow \lambda \circ P \circ \tau,$ we obtain a unique polynomial $P=T_6.$
$\square$

\medskip

Let $V, V_1 \subset {\cal P}$ be fixed linear spaces. Let us denote by $ND_{V,V_1}$ the set of
polynomials $P\in V$ non-definite with respect to $V_1$. In particular, for
$V={\cal P}_d,V_1 = {\cal P}$ we denote the corresponding set by $ND_d$. If $V_1$ is a line
spanned by a fixed $Q\in {\cal P}$ we write $ND_{V,V_1}$ as $ND_{V,Q}$.
\bp \label{dimnd}
For each $V_1 \subset {\cal P}$ and $V \subset{\cal P}_d$ we have $ND_{V,V_1}\subset ND_d.$
The dimension of $ND_d$ does not exceed $[{d\over 6}]+1$.
\ep
\pr
The inclusion is immediate: any polynomial non-definite with respect to a smaller subspace
is non-definite with respect to a larger one. By Theorem \ref{definite} the set $ND_d$
consists of all $P \in {\cal P}_d$ which have $s\geq 2$ mutually $[a,b]$-prime right
$[a,b]$-factors. Hence $ND_d \subset \cup_{s\geq 2} DEC^d_s(a,b).$ By Proposition \ref{propp}
we have $\dim ND_d \leq [{d\over 6}]+1$. This completes the proof. $\square$.

\subsection{Co-definite polynomials}

Let $[a,b]$ and a subspace $V \subset {\cal P}_{[a,b]}$ be given.

\bt \label{ncd}
A polynomial $Q$ is not $V$-co-definite if and only if there exists a
polynomial $P \in V$ (necessarily non-definite) with a complete collection of
$[a,b]$-indecomposable right factors $W_1,\dots, W_s, \ s\geq 2,$
such that:

\smallskip

1. The polynomial $Q$ can be represented as $Q=\sum_{j=1}^s S_j(W_j),$

\smallskip

2. No one of $W_1,\dots, W_s$  is a right $[a,b]$-factor of $Q$.
\et
\pr
By Definition \ref{def.codef} a polynomial $Q$ is not $V$-co-definite if and only if
there exists a polynomial $P\in V$ such that all the moments $m_k=\int_a^b P^k(x)q(x)dx,$
$k\geq 0,$ vanish while $P$ and $Q$ do not satisfy the Composition Condition. Clearly, if
such $P$ exists it cannot be definite. Furthermore, by Theorem \ref{yui} the polynomial
$Q$ can be represented as a sum $Q=\sum_{j=1}^s S_j(W_j).$ Finally, since $P$ and $Q$ do
not satisfy Composition condition no one of $W_1,\dots, W_s$ can be an $[a,b]$-right factor
of $Q$.

In the opposite direction, assume that $P \in V$ as required exists. Since $Q$ possesses
a representation $Q=\sum_{j=1}^s S_j(W_j)$, where $W_1,\dots, W_s$ are right $[a,b]$-factors
of $P,$ we conclude (by linearity of the moments in $Q$) that all the moments
$m_k$, $k\geq 0,$ vanish. Furthermore, since $W_1,\dots, W_s$ is a complete collection of
right $[a,b]$-factors of $P$, the second assumption implies that $P$ and $Q$ do
not satisfy the Composition Condition. Hence $Q$ is not $V$-co-definite. $\square$

\bd \label{sums} For $V \subset {\cal P}$ we define the set ${\cal S}_{V,d} \subset {\cal P}_d$
as the set of polynomials $Q \in {\cal P}_d$ which can be represented as
$Q=\sum_{j=1}^s S_j(W_j),$ where $W_1,\dots, W_s$ are all $[a,b]$-indecomposable right
factors of a certain non-definite $P\in V.$ The set ${\cal S}_V$ is the union
$\cup_d {\cal S}_{V,d}$.
\ed
By Theorem \ref{ncd}, in order to describe explicitly all $V$-co-definite polynomials
up to degree $d$ we have first to describe the set ${\cal S}_{V,d}$ and then to describe
those $Q \in {\cal S}_{V,d}$ for which no one of $W_1,\dots, W_s$ is a right
$[a,b]$-factor of $Q$. Both these questions in their general form turn out to be rather
tricky, and we provide here only very partial results.

let us stress that the role of Theorem \ref{noncodef} below and of the rest of the results in this section is to describe in detail
the set of {\it non-composition} solutions of the Moment equations.  This will prepare an application, in Section 
\ref{mainr} below, of the second set of the Center Equations at infinity, according to Theorem \ref{ceqinf}: those 
provided by the vanishing of the second Melnikov function. We expect that together the two sets of equations always
imply Composition Condition (compare Conjecture 2 in Section \ref{mainr} below.)

To make formulas easier, without loss of generality we shall assume that $[a,b]$ coincides with
$[-{{\sqrt 3}\over 2},{{\sqrt 3}\over 2}]$.

\bt \label{noncodef} Let $V={\cal P}_{9,[-{{\sqrt 3}\over 2},{{\sqrt 3}\over 2}]}$. Then
the set ${\cal S}_{V,d}$ is a vector space consisting of all polynomials $Q\in {\cal P}_d$
representable as $Q=S_1(T_2)+S_2(T_3)$ for some polynomials $S_1,$ $S_2$. Furthermore,
the dimension ${\cal S}_{V,d}$ is equal to
$[{{d+1}\over 2}]+[{{d+1}\over 3}]-[{{d+1}\over 6}].$ In particular, this dimension
does not exceed $[{2\over 3}d]+1.$

For $d\leq 4$ the space ${\cal S}_{V,d}$ coincides with ${\cal P}_d$, and starting
with $d=5$ this space is always a proper subset of ${\cal P}_d$. We have
${\cal S}_{V,5}={\cal P}_4 \subset {\cal P}_5$ and ${\cal S}_{V,6}$ is the subspace
in ${\cal P}_6$ consisting of all the polynomials $Q$ of the form $Q= Q_1+\alpha T_3$
with $Q_1$ even of degree at most $6$. ${\cal S}_{V,7}={\cal S}_{V,6},$ while ${\cal S}_{V,8}$
is the subspace in ${\cal P}_8$ consisting of all the polynomials $Q$ of the form
$Q= Q_1+\alpha T_3$ with $Q_1$ even of degree at most $8$. ${\cal S}_{V,9}$
is the subspace in ${\cal P}_9$ consisting of all the polynomials $Q$ of the form
$Q= Q_1+\alpha T_3+\beta T^3_3$ with $Q_1$ even of degree at most $8$.
\et
\pr
By Theorem \ref{cla} the only non-definite polynomials in
$V={\cal P}_{9,[-{{\sqrt 3}\over 2},{{\sqrt 3}\over 2}]}$ are scalar multiples of $T_6$.
$T_6=T_2\circ T_3=T_3\circ T_2$ has exactly two right
$[-{{\sqrt 3}\over 2},{{\sqrt 3}\over 2}]$-factors $T_2$ and $T_3$. This proves the first
claim of Theorem \ref{noncodef}.

Next observe that
\be \label{pere}
\mathbb C[T_n]\cap \mathbb C[T_m]= \mathbb \C[T_l],
\ee
where $l={\rm lcm}(n,m).$ Indeed if $P$ is contained in $\mathbb C[T_n]\cap \mathbb C[T_m]$,
then there exist polynomials $A,$ $B$ such that
$$
P=A\circ T_n=B\circ T_m,
$$
and in order to show that there exists a polynomial $U$ such that $P =U \circ T_l$ one can
use the second Ritt theorem. However, such a proof is more difficult that it seems since it
requires an analysis of the possibility provided by \eqref{sol1} (see e.g. Lemma 4.1 of \cite{pak}).
It is more convenient to observe  that identity  \eqref{fifth} implies that the function
$$
F=P\circ \frac{1}{2}\left(z+\frac{1}{z}\right)=A\circ\frac{1}{2}\left(z^n+\frac{1}{z^n}\right)
=B\circ \frac{1}{2}\left(z^m+\frac{1}{z^m}\right)
$$
is invariant with respect to the both groups $D_{2n}$ and $D_{2m}$, where $D_{2s}$ is the
dihedral group generated by the transformations $z\rightarrow 1/z$ and $z\rightarrow e^{2\pi i/s} z.$
Therefore, $F$ remains invariant with respect to the group $<D_{2n},D_{2m}>=D_{2l}$ generated by $D_{2n}$ and $D_{2m}$, 
implying that there exists a rational function $U$ such that
$$
F=U\circ \frac{1}{2}\left(z^l+\frac{1}{z^l}\right).
$$
Since
$$
U\circ \frac{1}{2}\left(z^l+\frac{1}{z^l}\right)=U\circ T_l\circ \frac{1}{2}\left(z+\frac{1}{z}\right),
$$
we conclude that $P=U\circ T_l$, and it is easy to see that $U$ actually is a polynomial.

Denote by $U_{d,n}$ the subspace of $\mathbb C[T_n]$ consisting of all polynomials
of degree $\leq d.$ By the remark above we have
$U_{d,n}\cap U_{d,m}=U_{d,l}$. This implies that
$$
\dim {\cal S}_{V,d}=\dim U_{d,2}\oplus U_{d,3}-2=\dim U_{d,2}+\dim U_{d,3}-\dim U_{d,6}-2=
$$
$$
=[{{d+1}\over 2}] + [{{d+1}\over 3}]-[{{d+1}\over 6}]-1\leq [{2\over 3}d]+1.
$$
To get an explicit description of ${\cal S}_{V,d}$ for $d\leq 9$ we shall use the following
simple lemma, which is used also in Section \ref{mainr}. Consider polynomials
$\hat T_2(x)=2x^2-{3\over 2}$, $\hat T_3(x)=T_3(x)=4x^3-3x$, and $\hat T_6=\hat T_2\circ \hat T_3$.
Our polynomials $\hat T_j, \ j=2,3,6,$ differ from the usual Chebyshev polynomials only in a constant term,
chosen in such a way that $\hat T_j$ vanish at the points $-{{\sqrt 3}\over 2},{{\sqrt 3}\over 2}.$
In the representation
\be \label{sum2}
Q=S_1(\hat T_2)+S_2(\hat T_3), \ S_1,S_2 \in {\cal P}.
\ee
$Q$ belongs to ${\cal P}_d$, so it vanishes at the points $-{{\sqrt 3}\over 2},{{\sqrt 3}\over 2}.$ Hence we
can assume that both $S_1$ and $S_2$ do not have constant terms. Next we notice that all the even polynomials
$Q$ in ${\cal P}_d$ and only them are representable as $Q=S_1(\hat T_2)$.

\smallskip

Let $\deg S_1=m, \ \deg S_2=n$.

\bl \label{degrees}
Let $Q \in {\cal P}_d$ be represented via (\ref{sum2}). Then the polynomials $S_1$ and $S_2$ in (\ref{sum2})
can be chosen in such a way that $S_2$ is odd, and $\max (2m,3n)\leq d$.
\el
\pr
It is enough to consider only odd polynomial $S_2$. Indeed, it is immediate that all the even polynomials,
and only them are representable as $S_1(\hat T_2)$. Since for $l$ even $\hat T_3^l=x^l(4x^2-3)^l$ is an even
polynomial (and it is odd for $l$ odd), all the even degrees in $S_2$ can be omitted.

\smallskip

Under this assumption, the odd degree $n$ of $S_2$ must satisfy $3n\leq d$. Indeed, otherwise the odd degree
$3n$ of $S_2(\hat T_3)$ would be larger than $d$, and the highest degree term in this polynomials could not
cancel with the terms of $S_1(\hat T_2)$. By the same reason, assuming that $S_2$ is odd, we conclude that
$2m\leq d$. Indeed, otherwise the even degree $2m$ of $S_1(\hat T_2)$ would be larger than $d$, and the highest
degree term in this polynomials could not cancel with the terms of $S_2(\hat T_2)$. $\square$

\smallskip

Application of Lemma \ref{degrees} completes the proof of Theorem \ref{noncodef}. $\square$

\bt A polynomial of the form $Q=S_1(T_2)+S_2( T_3),$ where $S_1,$ $S_2$ are non-zero
polynomials, has $T_2$ (resp. $T_3$) as its right factor if and only if $S_2$ is a polynomial
in $T_2$ (resp. $S_1$ is a polynomial in $T_3$).
\et
\pr Indeed, assume say that $S_1( T_2)+S_2(T_3)=R(T_2)$ for some polynomial $R$. Then by
\eqref{pere} there exists a polynomial $F$ such that
$$S_2\circ T_3=F\circ T_6=F\circ T_2\circ T_3$$ implying that $S_2=F\circ T_2$.
$\square$

\bc Let $V={\cal P}_{9,[-{{\sqrt 3}\over 2},{{\sqrt 3}\over 2}]}$. A polynomial
$Q\in {\cal P}_{8,[-{{\sqrt 3}\over 2},{{\sqrt 3}\over 2}]}$ is not $V$-co-definite
if and only if it can be represented in the form
\be \label{kot}
Q=R+\alpha T_3,\ \ \ \alpha \in \C,
\ee
where $\alpha\neq 0$, and $R\in {\cal P}_{8,[-{{\sqrt 3}\over 2},{{\sqrt 3}\over 2}]}$
is an even polynomial distinct from $\beta T_6+\gamma,$ $\beta, \gamma \in \C.$
\ec
\pr
By the above results, if $P\in {\cal P}_8$ is not co-definite it can be represented in
the form $Q=S_1(T_2)+S_2(T_3),$ where $\deg S_1 \leq 4,$ and $S_1$ is not a linear
polynomial in $ T_3$, while $\deg S_2 \leq 2,$ and $S_2$ is not a linear polynomial in
$T_2.$ Since $S_2$ can be represented in the form $\delta T_2+\alpha z +\kappa,$ where
$\delta, \alpha, \kappa \in \C$, we conclude that such $Q$ can be represented in the form
\be \label{kot1}
Q=\tilde S_1(T_2)+\alpha T_3,
\ee
where $\deg \tilde S_1 \leq 4$. Furthermore, $\alpha\neq 0$, since otherwise $Q$ is
a polynomial in $T_2$, and $\tilde S_1$ is not a linear polynomial in $T_3,$
since otherwise $Q$ is a polynomial in $T_3.$ Therefore, since $\C[T_2]=\C[z^2]$
and $T_3\in {\cal P}_8$, the polynomial $P$ admits representation \eqref{kot}.

In other direction, it follows from \eqref{kot} that  \eqref{kot1} holds, where
$\alpha\neq 0$ and $\tilde S_1\neq \beta T_3+\gamma,$ $\beta, \gamma \in \C,$
implying that $Q$ is not co-definite. $\square$

\subsection{Polynomials with a special structure}

Let ${\cal R}=\{r_1,r_2,\dots\}$ be a set of prime numbers, finite or infinite.
Define $U({\cal R})$ as a subset of ${\cal P}$ consisting of polynomials
$P=\sum_{i=0}^N a_i x^i$ such that for any non-zero coefficient $a_i$ the degree $i$
is either coprime with each $r_j\in {\cal R}$ or it is a power of some $r_j\in {\cal R}$.
Similarly, define $U_1({\cal R})$ as a subset of ${\cal P}$ consisting of polynomials $Q$
such that for any non-zero coefficient $a_i$ of $Q$ all prime factors of $i$ are contained
in ${\cal R}$. In particular, if ${\cal R}$ coincides with the set of all primes numbers,
then $U({\cal R})$ consists of polynomials in ${\cal P}$ whose degrees with non-zero
coefficients are powers of primes, while $U_1({\cal R})={\cal P}.$

\bt \label{speccase}
Let ${\cal R}=\{r_1,r_2,\dots\}$ be fixed. Then for any $a\ne b$ each polynomial
$P\in U({\cal R})$ is $[a,b]$-definite, and, in particular, it is $[a,b]$-definite
with respect to $U_1({\cal R})$, and each $Q\in U_1({\cal R})$ is $[a,b]$-co-definite
with respect to $U({\cal R})$.
\et
\pr
We show that vanishing of all the moments $m_k=\int^b_a P^k(x)q(x)dx$ for
$P\in U({\cal R})$ and $Q\in U_1({\cal R})$ implies Composition Condition.
By the construction, the degree of any $Q \in U_1({\cal R})$ is the
product of certain prime numbers in ${\cal R}$. By Corollary 4.3  of \cite{pp} vanishing
of the moments implies that the degrees of $P$ and $Q$ cannot be mutually prime. Hence
$\deg P$ is divisible by one of $r_j$. But then by the definition
this degree must be a power of $r_j$. Finally, it was shown in \cite{pp} (see also
Section 3.2.1 above) that polynomials $P$ with $\deg P$ a power of a prime number are
definite. Hence vanishing of the moments $m_k$ implies Composition condition for $P,Q$
on $[a,b]$. $\square$

\subsection{The Moment and the Composition Sets}

Using the information on definite and co-definite polynomials provided above we now
can describe more accurately the interrelation between the Moment and the Composition
sets.

Let $V, V_1 \subset {\cal P}$ be fixed linear spaces. As above, $ND_{V,V_1}$ is the set
of polynomials $P\in V$ non-definite with respect to $V_1$.

\bt \label{momdef} For each $Q \in V_1$ we have $MS_{V,Q}=COS_{V,Q} \ \cup \ N$ where $N$
is contained in $ND_{V,V_1} \subset ND$. In particular, for $V \subset {\cal P}_d$ and
any $Q$ the dimension of $N$ is at most $[{d\over 6}]+1.$
\et
\pr
If $P \in MS_{V,Q}$ but $P$ is not in $COS_{V,Q}$ then $P$ is not definite with respect
to $V_1$, and hence it belongs to $ND_{V,V_1},$ which is always a subset of $ND$.
If $V \subset {\cal P}_d$ then $P\in ND_d$ and the bound on the dimension follows from
Proposition \ref{dimnd}. $\square$

\smallskip

\noindent{\bf Example} (\cite{broy}) Let $[a,b]=[-{\sqrt 3\over 2},
{\sqrt 3\over 2}].$ Put $Q=(T_2+T_3)$, and consider $V = {\cal P}_6.$ Then the Moment set
$MS_{V,Q}$ contains exactly two components: the composition component
$COS_{V,Q} = \{P=R(T_2+T_3)\}$, with $R$ any polynomial of degree $2$, and the
non-composition component ${\cal T}=\{P=\a T_6, \ \a \in {\mathbb C} \}$. Here ${\cal T}$,
in fact, coincides with $ND_{V,Q}$.

\medskip

Our description of co-definite polynomials in Section 4.2 produces the following
result on the Moment and Composition Sets:
\bc \label{mscos1}
Let $V \subset {\cal P},$ and let $V_1 \subset {\cal P}_d$ be such that
$V_1\cap {\cal S}_{V,d}=\{0\}$, in the notation of Definition \ref{sums}.
Then for each $Q \in V_1$ we have $ND_{V,V_1}=\emptyset$ and $MS_{V,Q}=COS_{V,Q}$.
\ec
\pr
By Theorem \ref{ncd} and via Definition \ref{sums}, each $Q\in V_1, \ Q\ne 0$ is
co-definite with respect to $V$. Consequently, each $P\in V$ is definite with respect
to such $Q$. Application of Theorem \ref{momdef} completes the proof. $\square$

\medskip

In situation of Section 4.3 we get
\bc \label{mscos2}
For a fixed set $\cal R$ of prime numbers let $V=U({\cal R}), \ V_1=U_1({\cal R})$,
in notations of Section 4.3. Then for each $Q \in V_1$ we have $MS_{V,Q}=COS_{V,Q}$.
\ec
\pr
The result follows directly from Theorems \ref{momdef} and \ref{speccase}. $\square$

\section{Center Set near infinity}
\setcounter{equation}{0}

Let a polynomial $Q$ and a linear subspace $V \subset {\cal P}_d$ be fixed.
In this section we analyze the structure of the Center Set $CS_{V,Q}$ at and
near the infinite hyperplane $HV$, as compared to the Moment and Composition Sets
$MS_{V,Q}$ and $COS_{V,Q}$.  By Proposition \ref{inclusion} we have at infinity
$\bar {COS} \subset \bar {CS} \subset \bar {MS}.$

An important fact is that for each {\it definite} $P_0 \in \bar {CS}$ there is an entire
projective neighborhood $U$ of $P_0$ in $PV$ where $CS$ and $COS$ coincide:

\bt \label{csinf.def}
Let $P_0 \in \bar {CS}_{V,Q}$ be a definite polynomial. Then

\smallskip

1. In fact, $P_0 \in \bar {COS}_{V,Q}.$

\smallskip

2. There exists a projective neighborhood $U$ of $P_0$ in $P V$ such that
$CS_{V,Q} \cap U=COS_{V,Q} \cap U.$

\smallskip

3. $CS_{V,Q} \cap U$ is a linear space defined by vanishing of the linear parts
of the Center Equations. In particular, $CS$ is regular in $U$ and its local ideal is
generated by the Center Equations.
\et
\pr
From the inclusion $\bar {CS} \subset \bar {MS}$ we get $P_0 \in \bar {MS}_{V,Q}.$
Since the polynomial $P_0$ is definite by the assumptions, moments vanishing for this
polynomial implies composition, so in fact $P_0 \in \bar {COS}_{V,Q}.$

\smallskip

In homogeneous coordinates $(P,\nu)$ in $PV$ near $P_0$ put
$P=P_0+P_1, \ P_1\in V.$  By Proposition 7.2 of \cite{broy} the only nonzero linear terms
in the expansions of the homogenized Center Equations around the point $(0,0)$ in
variables $P_1, \ \nu$ are given by the following linear functionals in $P_1$:

\be
L_k(P_1) = -(k-3)\int^b_a P^{k-4}_0(x)q(x)P_1(x)dx, \ k=4,5,\dots .
\ee
Denote by $L \subset V$ the subspace defined by the linear equations
$L_k(P_1)=0, \ \ k=4,5,\dots$. Let us show first that $L\subset COS_{V,Q}$. Consider
certain polynomial $P_1 \in L$. Since $P_0$ is definite, vanishing of $L_k(P_1)$ implies
composition condition for $P_0(x)$ and $S(x)=\int^x_a P_1(\tau)q(\tau) d\tau$. Since,
being definite, $P_0$ has only one $[a,b]$-prime right composition $[a,b]$-factor $W$,
we conclude that $S=\tilde S(W)$. By the same reason, from $P_0\in COS_{V,Q}$ it follows
that $Q=\tilde P(W)$. Now Lemma 7.3 in \cite{broy} implies that $P_1=\tilde P_1(W)$, i.e.
$P_1\in COS_{V,Q}$, and hence $L\subset COS_{V,Q}$.

It follows that all the Center Equations vanish on $L$, which is the zero set of their
linear parts. Now we are in a situation of Lemma 7.4 of \cite{broy} (Nakayama Lemma in
Commutative algebra - see for example \cite{gg}, chapter 4, lemma 3.4). The conclusion
is that $CS=L=COS$ in a neighborhood of $P_0$, and the local ideal of this set is
generated by the Center Equations. This completes the proof of Theorem \ref{csinf.def}.
$\square$

\section{Main results} \label{mainr}
\setcounter{equation}{0}

Let $a\ne b$ be fixed. Below we denote by $\tilde T_j$ the transformed Cebyshev
polynomials $\tilde T_j=T_j\circ \mu,$ $\mu$ being a linear polynomial transforming
the couple $(a,b)$ to the couple $(-{{\sqrt 3}\over 2},{{\sqrt 3}\over 2})$.

Let linear subspaces $V, V_1 \subset {\cal P}_{[a,b]}$ and a polynomial $Q\in V_1$ be
fixed. The affine Center Set $CS_{V,Q}$ always contains the Composition Set $COS_{V,Q}$. In
this section we provide an upper bound for the dimensions of affine {\it non-composition}
components in $CS$. As above, $ND_{V,V_1} \subset ND$ denotes the set of $V_1$ non-definite
polynomials in $V$. For each affine algebraic set $A \subset V$ let $\bar A$ denote the
intersection of $A$ with the infinite hyperplane $HV$.

\bt \label{csinf1}
For each irreducible non-composition component $A$ of the affine Central Set
$CS_{V,Q}$ we have $\bar A \subset \bar {CS}_{V,Q}\cap ND \subset \bar {MS}_{V,Q} \cap ND$.
Consequently, $\dim A \leq \dim (\bar {MS}_{V,Q} \cap ND)+1$. In particular,
for any polynomial $Q$, and $V \subset {\cal P}_d$ the dimension of $A$ cannot exceed
$[{d\over 6}]+2$.
\et
\pr
We always have $\bar A \subset \bar CS \subset \bar MS.$ Now, if $P_0 \in \bar A$ then
$P_0$ cannot be definite. Indeed, otherwise there would exist a neighborhood $U$ of $P_0$
provided by Theorem \ref{csinf.def}, where $A \cap U \subset COS \cap U$. Since $A$
is irreducible, this would imply that $A \subset COS$, which contradicts the assumption
that $A$ is a non-composition component of $CS$. Thus
$\bar A \subset \bar {MS}_{V,p} \cap ND$. Now since the infinite hyperplane $HV$ has
codimension one in the projective space $PV$,
for each $A$ we have $\dim A \leq \dim \bar A+1$. Application of Proposition \ref{dimnd}
completes the proof. $\square$.

\medskip

Notice that the dimension of the composition components of $CS$ may be of order $d\over 2$,
while by Theorem \ref{csinf1} the dimension of the non-composition components is of order
at most $d\over 6.$ To our best knowledge, this is the first general bound of this form for
the polynomial Abel equation.

\bc \label{csinf2} (\cite{broy}).
Let $V ={\cal P}_5$. Then for any $Q$ the Center Set $CS_{V,Q}$ consists of a
Composition Set with possibly a finite number of additional points.
\ec
\pr
By Theorem \ref{cla} there are no non-definite polynomials in $V ={\cal P}_5$. So the set
$\bar {MS}_{V,Q} \cap ND$ is empty and its dimension is $-1$. $\square$.

\bc \label{csinf3}
Let $V ={\cal P}_9$. Then for any $Q$ the Center Set $CS_{V,Q}$ consists of a Composition
Set with possibly a finite number of additional curves.
\ec
\pr
By Theorem \ref{cla} the only non-definite polynomials in $V ={\cal P}_9$ are scalar
multiples of $\tilde T_6$. So the set $\bar {MS}_{V,Q} \cap ND$ consists at most of one
point, and its dimension is at most $0$. $\square$

\bc \label{csinf4}
Let $V ={\cal P}_{11}$. Then for any $Q$ the  $CS_{V,Q}$ consists of a
Composition Set with possibly a finite number of additional two-dimensional components.
\ec
\pr
Theorem \ref{cla} describes non-definite polynomials in $V ={\cal P}_{11}$.
We see that the set $\bar {MS}_{V,Q} \cap ND$ consists at most of a finite
number of points, and a one-dimensional component, and its dimension is at most $1$.
$\square$.

\medskip

Notice that the bounds of Corollaries \ref{csinf2}-\ref{csinf4} are more accurate than
the general bound of Theorem \ref{cla}.

\medskip

Recall that by Definition \ref{sums} the set ${\cal S}_V$ consists of all $Q$ which can
be represented as $Q=\sum_{j=1}^s S_j(W_j),$ where $W_1,\dots, W_s$ are all
$[a,b]$-indecomposable right factors of a certain $P\in V.$

\bt \label{csinf5}
Let $V \subset {\cal P},$ and let $Q \in {\cal P}\setminus {\cal S}_V.$ Then the
Center Set $CS_{V,Q}$ consists of a Composition Set with possibly a finite number of
additional points. In particular, this is true for $V={\cal P}_9$ and any $Q$ not
representable as $Q=S_1(\tilde T_2)+S_2(\tilde T_3)$.
\et
\pr
This result follows directly from Theorem \ref{cla} and Corollary \ref{mscos1}. The
case $V={\cal P}_9$ is covered by Theorem \ref{noncodef}. However, since Theorem
\ref{csinf5} is one of the central results of this paper, we give its short independent
proof. We show that the Moment Set $MS_{V,Q}$ does not contain non-definite polynomials.
Indeed, for each non-definite $P\in V$ vanishing of the moments $m_k=\int_a^b P^k(x)q(x)dx$
implies $Q\in {\cal S}_V,$ by Definition \ref{sums}. But by our assumptions
$Q \in {\cal P}\setminus {\cal S}_V.$ Therefore $P$ is not in $MS_{V,Q}$. Application
of Theorem \ref{csinf1} completes the proof. $\square$.

\medskip

We expect that the result of Theorem \ref{csinf5} can be significantly extended. In particular,
we expect that the following statement is true:

\medskip

\noindent{\bf Statement 6.1} {\it Let $V \subset {\cal P}.$ Assume that either
$Q\in {\cal P}\setminus {\cal S}_V,$ or $Q \in {\cal S}_V,$ and it is not $V$-co-definite.
Then the Center Set $CS_{V,Q}$ consists of a Composition Set with possibly a finite number
of additional points.}

\smallskip

Closely related to Statement 6.1 is the following

\medskip

\noindent{\bf Conjecture 2} {\it For polynomials $P,Q$ vanishing of all the moments
$m_k(P,Q)$ and of all the second Melnikov coefficients $D_j(P,Q)$ (see Theorem \ref{ceqinf})
implies Composition Condition.}

\bt \label{csinf6}
Conjecture 2 implies Statement 6.1.
\et
\pr
Assume, as in Statement 6.1, that either $Q\in {\cal P}\setminus {\cal S}_V,$ or
$Q \in {\cal S}_V,$ and it is not $V$-co-definite. The first case is treated in Theorem
\ref{csinf5}. In the second case we still show that the Center Set at infinity $\bar {CS}_{V,Q}$
does not contain non-definite polynomials. Assume, in contradiction, that $P\in \bar {CS}_{V,Q}$
is non-definite, and let $W_1,\dots, W_s, \ s\geq 2,$ be all the $[a,b]$-indecomposable right
factors of $P$. According to Theorem \ref{ceqinf} $P$ satisfies equations $m_k(P,Q)=0$ and
$D_j(P,Q)=0$. By the first set of these equations $Q=\sum_{j=1}^s S_j(W_j),$ and by the second
set and by Conjecture 2 we conclude that one of $W_j$ is a right factor of $Q$. Now according
to Theorem \ref{ncd} $Q$ is $V$-co-definite, in contradiction with the assumptions. This
completes the proof. $\square$

\smallskip

Our next result confirms Conjecture 2, and hence Statement 6.1 for $\deg P, \deg Q \leq 9.$

\bt \label{csinf7}
1. Conjecture 2 is valid for $\deg P, \deg Q \leq 9,$ i.e. vanishing of all the moments $m_k(P,Q)$ and of
four initial second Melnikov coefficients $D_j(P,Q)$ implies Composition Condition for such $P,Q$.

\smallskip

2. Consequently, for $V ={\cal P}_9$, and for any $Q$ of degree up to $9$ not of the form
$Q=S_1(\tilde T_2)+S_2(\tilde T_3)$, or of this form, but such that neither $\tilde T_2$ nor
$\tilde T_3$ are the right composition factors of $Q$, the center set $CS_{V,Q}$ consists of
a composition set with possibly a finite number of additional points.
\et
\pr
By Theorem \ref{csinf6}, the first part of Theorem \ref{csinf7} implies its second part. So let polynomials $P,Q$ with
$\deg P, \deg Q \leq 9$ be given. If $P\ne \alpha \tilde T_6,$ it is definite, and hence already vanishing of
all the moments $m_k(P,Q)$ implies Composition Condition for $P,Q$. Consider now the case $P=\tilde T_6$. Here vanishing
of $m_k(P,Q)$ implies that $Q$ has a form $Q=S_1(\tilde T_2)+S_2(\tilde T_3)$ for some polynomials $S_2$ and $S_3$. By
Lemma \ref{degrees} we conclude that $S_1,S_2$ can be written in the form
$S_1(T)=\sum_{i=1}^4 c_i T^i, \ S_2(T)=\sum_{i=1}^2 \alpha_i T^{2i-1}$. Now we use the second set of the Center equations
at infinity: $D_j(P,Q)=0$.
\bp \label{meln1} The first four equations at infinity $D_j(P,Q)=0$ given in Theorem \ref{ceqinf} can be written as
\begin{eqnarray}\label{equat1}
\begin{split}
&D_1(P,Q)=\int_a^b Q^2p&=0,\\
&D_2(P,Q)=\int_a^b Q^2Pp&=0,\\
&D_3(P,Q)=2\int_a^b Q^2P^2p+\int_a^b Q(t)P(t)p(t)dt\int_a^t Q(\tau)p(\tau)d\tau&=0.\\
&D_4(P,Q)=\int_a^b Q^2P^3p+\int_a^b Q(t)P^2(t)p(t)dt\int_a^t Q(\tau)p(\tau)d\tau&=0.\\
\end{split}
\end{eqnarray}
\ep
\pr The proof is based on rather lengthy computations. We shall use the following result from \cite{broy}:

\bt\label{2.2.7} (Thorem 2.2, \cite{broy}) Any iterated integral $I_\a$ with $m_0+m_1+m_2$ appearances of $p$ and exactly
two appearances of $q$ after $m_0$ and $m_1$ appearances of $p$, respectively, can be transformed via integration by parts
to the sum of the iterated integrals of the following form:

\be\label{iter.2}
I_\a=\sum_{i=0}^{m_1} {{(-1)^{m_0+m_1-i}}\over {m_0!m_2!i!(m_1-i)!}} \int_a^b P^{m_0+i}(x)q(x)dx \int_a^x P^{m_1+m_2-i}(t)q(t)dt.
\ee
\et
Below we present calculations of the Melnikov coefficients at infinity $D_1(P,Q)$ and $D_3(P,Q)$. $D_2$ and $D_4$ are obtained in
a similar way. Let us recall that 

\be\label{PQ}
P(a)=P(b)=Q(a)=Q(b)=0,
\ee
while from Theorem \ref{ceqinf} for $k=4,6,8,10$ it follows that

\be\label{PQq}
\int_a^b P^iq=0, \ i=1,...,4.
\ee
Case 1. Let $k=5$. Then the corresponding Center Equation at infinity is given by $D_1=\sum n_\a I_\a =0,$ where the sum runs
over all the indices $\a=(\a_1,...,\a_s)$ with exactly two appearances of $1$, and with $\sum_{j=1}^s \a_j=k-1=4$. Hence we have
exactly one appearance of $2$, and $s=3$.

Now, for $\a_1=(1,1,2)$ we have $n_{\a_1}=-12, \ m_0=m_1=0,m_2=1,$ and hence

$$
I_{\a_1}=\int_a^b q(x_1)dx_1\int_a^{x_1}q(x_2)dx_2\int_a^{x_2} p(x_3)dx_3 = 
$$ 
$$
\left [ Q(x_1)\cdot \int_a^{x_1}q \int_a^{x_2} p \right ] \ |_{x_1=a}^b - \int_a^b Qq dx_1 \int_a^{x_1}p = 
$$
$$
-{1\over 2} \left [ Q^2(x_1)\cdot \int_a^{x_1}p \right ] \ |_{x_1=a}^b +{1\over 2}\int_a^b Q^2p= {1\over 2}\int_a^b Q^2p.
$$
For $\a_2=(1,2,1)$ we have $n_{\a_2}=-8, \ m_0=0, m_1=1,m_2=0,$ and 
$$
I_{\a_2}=\int_a^b q\int_a^{x_1}p\int_a^{x_2} q = \left [ Q(x_1)\cdot \int_a^{x_1}p \int_a^{x_2} q \right ] 
\ |_{x_1=a}^b - \int_a^b Qp \int_a^{x_1}q =
$$ 
$$
- \int_a^b Q(x_1)p(x_1)(Q(x_1)-Q(a))dx_1=-\int_a^b Q^2p.  
$$
For $\a_3=(2,1,1)$ we have $n_{\a_3}=-6, \ m_0=1, m_1=0,m_2=0,$ and 
$$
I_{\a_3}=\int_a^b p\int_a^{x_1}q\int_a^{x_2} q = \int_a^b p\int_a^{x_1}Qq=
$$
$$
{1\over 2}\int_a^b (Q^2(x_1)-Q^2(a))p(x_1)dx_1={1\over 2}\int_a^b Q^2p.
$$
Then $D_1=\sum_{i=1}^3 n_{\a_i}I_{\a_i}=(-12\cdot {1\over 2} +8-6\cdot {1\over 2})\int_a^b Q^2p=\int_a^b Q^2p$.

\medskip

\noindent Case 2. For $k=7,9,11$ we shall use expressions (\ref{iter.2})-(\ref{PQq}), which will allow us to present somewhat 
lengthy calculations in a more systematic way. For $k=7$ the calculations easily provide the answer $D_2=\int_a^b Q^2Pp=0$, so 
we concentrate on the case $k=9$. The Center Equation at infinity is given by  
\be\label{k9}
\sum n_\a I_\a=0 
\ee
with $\sum_{j=1}^s \a_j=k-1=8$, and with exactly two appearances of $1$. Hence we have exactly three appearances of $2,$ and
$s=5$.

For $\a_1=(1,1,2,2,2)$ we have $n_{\a_1}=-8\cdot7\cdot5\cdot3, \ m_0=m_1=0, \ m_2=3$. Then by (\ref{iter.2}) we have 
$I_{\a_1}={1\over {3!}}\int_a^b q(x)dx \int_a^x P^3q.$ Integrating by parts once more we get $I_{\a_1}={1\over 4}\int_a^b Q^2P^2p$.

In a similar way, expressions (\ref{iter.2})-(\ref{PQq}) and some additional integrations by part allow us to express all $I_\a$
through $J_1=\int_a^b Q^2P^2p$ and $J_2=\int_a^b QPp\int_a^x Qp.$

\medskip

\noindent{\it For $\a_2=(1,2,1,2,2)$:} $n_{\a_2}=-8\cdot6\cdot5\cdot3, \ m_0=0, m_1=2, m_2=2$ and

$$
I_{\a_2}=-{1\over 2}\int_a^b q \int_a^x P^3q+{1\over 2}\int_a^b Pq \int_a^x P^2q=...=-{1\over 2}J_1-J_2.
$$
\noindent{\it For $\a_3=(1,2,2,1,2)$:} $n_{\a_3}=-8\cdot6\cdot4\cdot3, \ m_0=0, m_1=2, m_2=1$ and

$$
I_{\a_3}={1\over 2}\int_a^b q \int_a^x P^3q -\int_a^b Pq \int_a^x P^2q +{1\over 2}\int_a^b P^2q \int_a^x Pq=...=3J_2.
$$
\noindent{\it For $\a_4=(1,2,2,2,1)$:} $n_{\a_4}=-8\cdot6\cdot4\cdot2, \ m_0=0, m_1=3, m_2=0$ and

$$
I_{\a_4}=-{1\over 6}\int_a^b q \int_a^x P^3q +{1\over 2}\int_a^b Pq \int_a^x P^2q -{1\over 2}\int_a^b P^2q \int_a^x Pq+
$$
$$
{1\over 6}\int_a^b P^3q \int_a^x q=...=-2J_2.
$$
In the same way, for the remaining $6$ transpositions $\a_5=(2,1,1,2,2),\a_6=(2,1,2,1,2),\a_7=(2,1,2,2,1),\a_8=(2,2,1,1,2),
\a_9=(2,2,1,2,1),\a_{10}=(2,2,2,1,1)$ we obtain the following values of $(n_\a,I_\a)$:
$$
(-7\cdot6\cdot5\cdot3, -{1\over 4}J_1+J_2), \ (-7\cdot6\cdot4\cdot3, J_1-4J_2), \ (-7\cdot6\cdot4\cdot2, 3J_2), \
$$
$$
(-7\cdot5\cdot4\cdot3, -{1\over 4}J_1+J_2), \ (-7\cdot5\cdot4\cdot2, -{1\over 2}J_1-J_2), \ (-7\cdot5\cdot3\cdot2, {1\over 4}J_1).
$$
Substituting these expressions for $m_\a$ and $I_\a$ into (\ref{k9}) we finally obtain $-2(2J_1+J_2)$, so, after omitting a nonzero
coefficient $-2$, we get $D_3=2J_1+J_2=2\int Q^2P^2p+\int QPp\int Qp$.

The last equation in (\ref{equat1}), for $k=11$, is obtained in a completely similar way. $\square$

\medskip

The following results describe the application of these four equations to the specific combinations of Chebyshev
polynomials representing $Q$. To simplify the numeric coefficients we assume here that $[a,b]=[0,1]$ and so
$\tilde T_2(x)=x(x-1), \tilde T_3(x)=x(x-1)(2x-1)$. We also notice that $\tilde T_6=\tilde T_3^2=\tilde T_2^2+4\tilde T_2^3.$

\smallskip

In all the calculations below $P=\tilde T_6$ is fixed, while $Q=S_1(\tilde T_2)+S_2(\tilde T_3)$ with
$S_1(T)=\sum_{i=1}^4 c_i T^i, \ S_2(T)=\sum_{i=1}^2 \alpha_i T^{2i-1},$ is variable. We substitute these $P$ and $Q$ into the
equations (\ref{equat1}) of Proposition \ref{meln1} and get a system of algebraic equations with respect to the complex unknowns
$c_1,c_2,c_3,c_4,\alpha_1,\alpha_2$.

\smallskip

It is convenient to introduce the expressions $L_k=\int_0^1 S_1(\tilde T_2)\tilde T^k_3 d\tilde T_6$, which are linear forms
in $c_1,c_2,c_3,c_4$. Using these expressions we can rewrite equations (\ref{equat1}) as
\begin{eqnarray}\label{equat2}
\begin{split}
&\alpha_1L_1+\alpha_2L_3&=0,\\
&\alpha_1L_3+\alpha_2L_5&=0,\\
&{{16}\over {15}}\alpha_1L_5+{{36}\over {35}}L_7&=0,\\
&{{25}\over {21}}\alpha_1L_7+{{49}\over {45}}L_9&=0.\\
\end{split}
\end{eqnarray}

\bp\label{new.6.2}
The expressions $L_k, \ k=1,3,5,7,9$ can be written explicitly as the following linear forms in the coefficients $c_1,c_2,c_3,c_4$
of the polynomial $S_1(T)$:

$$
L_1= -{8\over {13}}\cdot {{(5!)^2}\over {11!}} (-13c_1+4c_2-c_3+{4\over {17}}c_4)
$$
$$
L_3= - {3\over {14\cdot 9}}\cdot {{(8!)^2}\over {17!}}(-{{38}\over 3}c_1+4c_2-c_3+{{16}\over {69}}c_4)
$$
$$
L_5= - {4\over {33\cdot 25}}\cdot {{(11!)^2}\over {23!}}(-{{25}\over 2}c_1+4c_2-c_3+{{20}\over {87}}c_4)
$$
$$
L_7= - {{10}\over {11\cdot 13\cdot 31}}\cdot {{(14!)^2}\over {29!}}(-{{62}\over 5}c_1+4c_2-c_3+{{8}\over {35}}c_4)
$$
$$
L_9= - {{9}\over {13\cdot 17\cdot 37}}\cdot {{(17!)^2}\over {35!}}(-{{37}\over 3}c_1+4c_2-c_3+{{28}\over {123}}c_4)
$$
\ep
\pr
Straightforward computation of $L_k$ using the identities $\tilde T_6=\tilde T_3^2=\tilde T_2^2+4\tilde T_2^3, \ 
d\tilde T_6=2\tilde T_3 d\tilde T_3=2(\tilde T_2+6\tilde T_2^2)d\tilde T_2,$ and 
$\int_0^1 \tilde T_2^n(x)dx=(-1)^n\cdot {{(n!)^2}\over {(2n+1)!}}.$ $\square$

\medskip

Now we come back to system \eqref{equat2}. Let us start with the special case where $\alpha_2=0$.

\bp \label{meln3}
Let $P=\tilde T_6, Q=S_1(\tilde T_2) + \alpha_1 \tilde T_3$, with $S_1(T)=\sum_{i=1}^4 c_i T^i$. If the first three equations
(\ref{equat1}) of Proposition \ref{meln1} are satisfied, then either $Q=S_1(\tilde T_2)$ or $Q=c_2 \tilde T_6+\alpha_1 \tilde T_3$.
In each of these cases $Q$ has either $\tilde T_2$ or $\tilde T_3$ as a right composition factor.
\ep
\pr
Substitution to the equations (\ref{equat2}) gives the following system of equations on the coefficients
$\alpha_1, c_1,c_2,c_3,c_4$:
\begin{eqnarray}
\begin{split}
&\alpha_1(-13c_1+4c_2-c_3+{4\over {17}}c_4)&=0,\\
&\alpha_1(-{{38}\over 3}c_1+4c_2-c_3+{{16}\over {69}}c_4)&=0,\\
&\alpha_1(-{{25}\over {2}}c_1+4c_2-c_3+{{20}\over {87}}c_4)&=0.
\end{split}
\end{eqnarray}
The result follows immediately from this system. $\square$

\smallskip

Let us consider now the remaining case, where $\alpha_2\ne 0.$
\bp \label{meln4}
Let $P=\tilde T_6, Q=S_1(\tilde T_2) + \alpha_1 \tilde T_3+\alpha_2 \tilde T^3_3$, with $S_1(T)=\sum_{i=1}^4 c_i T^i$, and
$\alpha_2\ne 0.$ If all the four equations (\ref{equat1}) of Proposition \ref{meln1} are satisfied, then
$Q=c_2 \tilde T_6+\alpha_1 \tilde T_3+\alpha_2 \tilde T^3_3$, and hence $Q$ has $\tilde T_3$ as a right composition factor.
\ep
\pr
Substitution to the equations (\ref{equat2}) gives a system of equations on the coefficients
$\alpha_1, \alpha_2, c_1,c_2,c_3,c_4$ which, putting $K:= {{\alpha_1}\over {\alpha_2}}$, can be brought to the following form:

\begin{eqnarray}\label{equat3}
\begin{split}
&(-4199K-19)c_1+(323K+{3\over 2})(4c_2-c_3)+(76K+{8\over {23}})c_4&=0\\
&(-874K-5)c_1+(69K+{2\over 5})(4c_2-c_3)+(16K+{8\over {87}})c_4&=0\\
&(-40600K-252)c_1+(3248K+{{630}\over {31}})(4c_2-c_3)+({{2240}\over 3}K+{{144}\over {31}})c_4&=0\\
&(-7750K-49)c_1+(625K+{{147}\over {37}})(4c_2-c_3)+({{1000}\over 7}K+{{1372}\over {1517}})c_4&=0
\end{split}
\end{eqnarray}
System \eqref{equat3} contains four equations with respect to four variables $K,c_1,c_4$ and $t=4c_2-c_3$. It presents a
system of four linear equations with respect to $c_1,c_4,t$, but $K$ enters the coefficients. We shall show that for each
$K$ this system has only the trivial solution $c_1=c_4=t=0$. Indeed, existence of a non-trivial solution would imply
simultaneous vanishing of, for example, the determinants

$$
\Delta_1(K)={{21280}\over 3}K^3+ {{2736}\over {31}}K^2+{{1368}\over {4495}}K+{{24}\over {103385}}
$$
and
$$
\Delta_2(K)= 76000K^3+ {{101998240}\over {104673}}K^2+{{3934112}\over {1081621}}K+{{3528}\over {1081621}},
$$
formed by the coefficients of system (6.3) in the rows $1,2,3$ and $1,3,4,$ respectively. However, the resultant
$res(\Delta_1(K),\Delta_2(K))$ of the polynomials $\Delta_1(K),\Delta_2(K)$ in $K$ is approximately equal to $21.51447438,$
so it is non-zero, and hence these polynomials do not have common roots. (The final calculation of $\Delta_1(K),\Delta_2(K)$
and of their resultant has been performed with the help of ``MATLAB'' system).

\smallskip

We conclude that for any $\alpha_1,\alpha_2$ with $\alpha_2\ne 0$, and for $K={{\alpha_1}\over {\alpha_2}}$, system (6.3)
implies $c_1=c_4=0, \ c_3=4c_2$. Consequently, any polynomial $Q$ satisfying this system has a form
$$
Q=c_2\tilde T^2_2+4c_2\tilde T^3_2+\alpha_1 \tilde T_3+\alpha_2 \tilde T^3_3=c_2\tilde T_6+\alpha_1 \tilde T_3+\alpha_2 \tilde T^3_3.
$$
In particular, $Q$ has $\tilde T_3$ as its right composition factor. This completes the proof of Theorem \ref{csinf7}: vanishing of
the moments and of the initial four Melnikov coefficients implies composition for $P,Q$ up to degree $9$. $\square$

\medskip

Finally we consider Center Sets in the subspaces $V=U_{\cal R}$, as defined in
Section 3.3.

\bt \label{csinf8}
Let a subset ${\cal R}=\{r_1,r_2,\dots\}$ of prime numbers be fixed. Put $V= U({\cal R})$,
as defined in Section 4.3 above. Then for any $a\ne b$ and for each fixed polynomial
$Q\in U_1({\cal R})$ the center set $CS_{V,Q}$ of Abel equation (\ref{Abel}) inside the space
$V$ consists of a Composition Set with possibly a finite number of additional points.
\et
\pr
This is a direct consequence of Corollary \ref{mscos2} and Theorem \ref{csinf1}. $\square$

\medskip

The results of this section cover all the results of Theorems \ref{Int1} - \ref{Int4}
stated in the Introduction.

\medskip

The methods developed in this paper work not only in the setting of the Center
Equations at infinity. They can be applied also to the study of the local structure
of the affine Center Set, extending the approach of \cite{bby}. Here we use ``second degree'' Nakayama lemma
in order to conclude that the Center Set (locally near the origin) coincides with the Composition Set, defined
by the Moments and the second Melnikov function. We plan to present these results separately. 

Our recent paper \cite{bpy} applies the results of the present paper on definite polynomials to the parametric 
versions of the Center-Focus problem for polynomial Abel equation. 

\vskip1cm

\bibliographystyle{amsplain}

\end{document}